%% file: TIsto_NLopt_Y_Ito_main_20250121d_submit_arxiv2.tex
\documentclass[]{article}  
\usepackage[margin=30truemm]{geometry}

\usepackage[cmex10]{amsmath} 
\usepackage{amssymb} 
\usepackage{amsthm}

\usepackage{cases} 

\usepackage{fixltx2e} 

\usepackage[subrefformat=parens]{subcaption} 
\usepackage{graphicx}
\usepackage{epstopdf}

\usepackage{xcolor}

\usepackage{cite}

\usepackage{algorithm}
\usepackage{algpseudocode}

\newtheorem{definition}{Definition}
\newtheorem{theorem}{Theorem}
\newtheorem{lemma}{Lemma}
\newtheorem{assumption}{Assumption}
\newtheorem{remark}{Remark}

\newcommand{\argmin}{\operatornamewithlimits{\mathrm{arg~min}}}

\definecolor{myGreen}{rgb}{ 0, 0.7, 0.3 }
\definecolor{myBlue}{rgb}{ 0, 0.4, 1 }
\definecolor{myPurple}{rgb}{ 0.7, 0, 0.3 }
\definecolor{myGray}{gray}{ 0.7 }
\definecolor{myFigGray}{gray}{ 0.5 }

\newcommand{\MyHighlight}[1]{\textbf{#1}}

\usepackage[overlay]{textpos}
\usepackage{pgfplots}
\pgfplotsset{compat=newest}
\usepackage{tikz}
\usetikzlibrary{positioning}
\usetikzlibrary{arrows}

\usepackage{natbib}

\begin{document}

	\title{Stabilizing Optimal Control for Nonlinear Stochastic Systems: A Parametric Gradient-Based Approach%
		\thanks{%
			This paper is submitted to a journal for possible publication. The copyright of this paper may be transferred without notice, after which this version may no longer be accessible.
			This work was partly supported by JSPS KAKENHI Grant Number JP18K04222.
			The materials of this paper have been published in part at \citep{ItoACC19}.
			We would like to thank Editage (www.editage.jp) for the English language editing.
	}} 
	
	\author{Yuji Ito\thanks{Yuji Ito is the corresponding author and with Toyota Central R\&D Labs., Inc., 41-1 Yokomichi, Nagakute-shi, Aichi 480-1192, Japan	(e-mail: ito-yuji@mosk.tytlabs.co.jp).}
		\and
		Kenji Fujimoto\thanks{Kenji Fujimoto is with the Department of Aeronautics and Astronautics, Graduate School of Engineering, Kyoto University, Kyotodaigakukatsura, Nishikyo-ku, Kyoto-shi, Kyoto 615-8540, Japan (e-mail: k.fujimoto@ieee.org).}
	}  
	
	\date{}
	
	\maketitle

	\begin{abstract} 
		This study proposes a method for designing stabilizing suboptimal controllers for nonlinear stochastic systems.
		These systems include time-invariant stochastic parameters that represent uncertainty of dynamics, posing two key difficulties in optimal control.
		Firstly, the time-invariant stochastic nature violates the principle of optimality and Hamilton-Jacobi equations, which are fundamental tools for solving optimal control problems.  
		Secondly, nonlinear systems must be robustly stabilized against these stochastic parameters.
		To overcome these difficulties simultaneously, this study presents a parametric-gradient-based method with a penalty function.
		A controller and cost function are parameterized using basis functions, and a gradient method is employed to optimize the controller by minimizing the parameterized cost function. 
		Crucial challenges in this approach are parameterizing the cost function appropriately and deriving the gradient of the cost.
		This study provides explicit formulations of an optimally parameterized cost and its gradient.
		Furthermore, a suitable penalty function is proposed to ensure robust stability, even when using the gradient method. 
		Consequently, the gradient method produces a suboptimal feedback controller that guarantees the robust stability.
		The effectiveness of the proposed method is demonstrated through numerical simulations, highlighting its performance in comparison with other baseline methods.
	\end{abstract}

\newcommand{\MyProof}{proof}

\input{TIsto_NLopt_Y_Ito_symbols_20241223.tex}

\section{Introduction}\label{sec_intro}

Various types of dynamical systems involve uncertainty and nonlinearity in the real world.
Uncertain physical parameters and/or noise in nonlinear dynamics are treated as follows: 
vehicular control systems such as
trajectory optimization with uncertain vehicle dynamics \citep{Listov24},
robot manipulators interacting in time-varying uncertain environments \citep{Liu24Cyb},
path planning for robotic spacecraft \citep{Nakka23},
sailboat navigation methods \citep{Miles22},
formation control for quadrotors \citep{Zhao22Cyb},
vehicle platooning \citep{Yin22}, and
robust suspension systems \citep{Bai22Cyb}.
Other applications include the following:
greenhouse production systems \citep{Svensen24},
chemical processes using continuous stirred-tank reactors \citep{Wu22},
mitigating pandemic waves \citep{Scarabaggio22},
individual differences in human motor systems using wearable devices \citep{Diaz23Cyb},
and nanoscale devices with manufacturing variations during mass production \citep{ItoACCESS19}.
Such uncertainty degrades the control performance and stability of systems.
These examples motivate us to consider stability and control of uncertain nonlinear systems, while the identification of the systems has also been addressed in existing work such as \citep{Ito23IFAC}.

Although the uncertainty of a system can be successfully represented by \textit{time-invariant} stochastic parameters, these parameters make nonlinear optimal control problems difficult.
The time-invariant parameters take random values obeying a probability distribution but are constant in time.
Because the time-invariant parameters violate the Markov property of system dynamics, they are not compatible with powerful tools for solving nonlinear optimal control, such as 
Hamilton-Jacobi-Bellman (HJB) equations \citep{Liu22Cyb}, 
the principle of optimality \citep{Lewis12},
generalized HJB equations \citep{Yang22Cyb}, and
Hamilton-Jacobi-Isaacs equations \citep{Fan22Cyb,Yang21Cyb}.
The details of this limitation are discussed in Remark \ref{rem:difficulties} in Section \ref{sec_problem_soc}.
Indeed, many stochastic optimal control methods, such as \citep{Ito24Automatica,ItoTAC19,Nakka23,Archibald20,McAllister23,Wu22}, have considered time-varying randomness that changes in time rather than time-invariant parameters.
If stochastic parameters are restricted to having discrete values in a finite set, Markov jump systems \citep{Chen22Cyb,Tao23Cyb} are elegant representations.
While generalized stochastic parameters have been considered \citep{Ito24TAC}, extending them to nonlinear systems remains future work.

Although various methods have tackled optimal control problems for time-invariant stochastic nonlinear systems, it is desirable but challenging to simultaneously satisfy crucial requirements.
The requirements include
guaranteeing \textit{robust stability} of the time-invariant stochastic nonlinear systems, 
developing \textit{feedback} controllers more robust than feedforward ones,
and designing the controllers \textit{offline} for real-time implementation.
Control design methods based on polynomial chaos (PC) \citep{Fisher09,Templeton09,Blanchard10,Ehlers11,Bhattacharya14}, stochastic model predictive control (MPC) \citep{Chen24}, and  gradient descent \citep{ItoCDC16,Ito24Cyb} have been applied only to linear systems.
The PC \citep{Xiu02} is efficient in treating stochastic uncertainty even for nonlinear control problems \citep{Listov24,Lefebvre20}.
Other techniques for nonlinear systems have been developed \citep{Phelps16a,Phelps16b,Okura16,Cottrill12}.
However, these nonlinear controllers  are feedforward rather than feedback.
Stochastic MPC \citep{Fagiano12,Lucia13,Mesbah14} is promising for realizing nonlinear feedback control.
Many approaches have been developed to implement stochastic MPC efficiently, such as 
linearization techniques \citep{Svensen24}, 
combining parameter estimation with robust adaptive MPC \citep{Sasfi23,Sinha22}, 
combining MPC with PC and Gaussian processes \citep{Bradford21}, 
Monte Carlo sampling for handling chance constraints \citep{Scarabaggio22},
a first-order second moment method to propagate uncertainty \citep{Yin22}, and 
particle MPC using many sample paths \citep{Dyro21}. 
Unfortunately, these MPC methods incur huge computational costs for online controller design, which restricts control applications.
Although recurrent neural networks have learned stochastic optimal controllers with low computational costs \citep{Yang22JPC}, they do not ensure optimality, feasibility, and closed-loop stability.
In other strategies, 
an existing method \citep{Pereira17b} has assumed that uncertain stochastic parameters are known or estimated.
A method with differential dynamic programming is limited to a region near predefined state trajectories \citep{Boutselis16}. 
Most importantly, the aforementioned control methods
have not guaranteed stability of nonlinear stochastic systems, indicating the difficulty of simultaneously realizing performance optimality and stabilization.

This study presents a method for solving optimal control problems for nonlinear systems containing \textit{time-invariant stochastic} parameters while simultaneously satisfying the aforementioned three requirements.
The proposed method provides \textit{feedback} controllers \textit{offline} that guarantee \textit{robust stability}.
Our strategy focuses on a gradient-based approach to minimize the expectation of a performance index regarding the stochastic parameters.
Although a simple gradient method is used to optimize a parameterized controller, it is challenging to obtain the gradient of the expected performance regarding the controller parameters. 
We derive this gradient to complete the controller design.
The proposed gradient method employs a penalty function based on sum of squares (SOS) to guarantee the robust stability of the nonlinear stochastic systems.
Our main contributions are listed below:

\begin{enumerate}
	\item
	\textbf{Simple and versatile concept:}
	To solve nonlinear optimal control problems involving time-invariant stochastic parameters,
	we propose a gradient-based approach with parametrization of a cost function and controller.
	Furthermore, the approach is combined with SOS to ensure robust stability of the systems regarding  the unknown stochastic parameters.
	The proposed approach provides a suboptimal controller satisfying the three requirements: feedback control, offline design, and robust stability, whereas the existing design methods have not realized them simultaneously.

	\item
	\textbf{Theory for parametric gradients:}
	This study designs an appropriate approximate cost function whose gradient is employed in the gradient method instead of a true gradient.
	Because a cost function over an infinite horizon is infeasible to evaluate in general, the infinite horizon cost is approximated as a parametric function such that a residual regarding optimality conditions is minimized.
	The relation between the residual and the approximate cost is clarified in Theorem \ref{thm:BellmanResidual}.
	Theorems \ref{thm:explicit_cost} and \ref{thm:explicit_grad} derive an optimal approximate cost function and its gradient in explicit form, respectively.
	The parametric function reduces to a linear combination of continuous basis functions such as kernel functions and polynomials.
	The resulting gradient is easy to implement because it is given as an explicit function of controller parameters.

	\item
	\textbf{Theory for robust stability:}
	This study guarantees the robust stability of nonlinear polynomial systems even with unknown stochastic parameters.
	We propose a penalty function used in the gradient method so that  an SOS-based stability analysis is combined with the minimization of the cost function.
	Theorem \ref{thm:stabilizing_optim} guarantees that the robust stability holds as long as the controller is successfully updated via the gradient method.
	Theorem \ref{thm:stabilizing_initial} derives an appropriate initial controller to start the gradient method.

	\item 
	\textbf{Demonstration:}
	The efficacy of the proposed approach is demonstrated using numerical examples.
	Average control performance is improved compared with three types of baseline controllers.
	The results support that the proposed method is successfully applied to nonlinear stochastic systems.

\end{enumerate}

This paper is an extension of the authors' conference proceedings \citep{ItoACC19} with the following main extensions:
The method used in this study guarantees robust stability, which cannot be guaranteed using \citep{ItoACC19}.  
A residual of approximating cost functions is clarified, whereas \citep{ItoACC19} has focused on the case of no residual.
Novel numerical examples are demonstrated to highlight the efficacy of the proposed approach.
Technical points have been modified so that the readability and technical soundness of this study are enhanced.

The remainder of this paper is organized as follows:
The mathematical notation is described in Section \ref{sec_notations}. 
Sections \ref{sec_problem} and \ref{sec_method} present the main problem and its solution, respectively.
Section \ref{sec_numerical_example} presents numerical examples that demonstrate the efficacy of the proposed approach.
Finally, Section \ref{sec_conclusion} concludes this study and discusses future work.

\section{Notation}\label{sec_notations}

The following notations are utilized in this paper:

\begin{itemize}

	\item 
	$\Identity{\DimANotation}$: the $\DimANotation \times \DimANotation$ identity matrix
	
	\item
	$\El{\NotationVec}{\IDNotation}$: the $\IDNotation$th component of $\NotationVec \in \mathbb{R}^{\DimANotation}$
	\item
	$\El{\NotationMat}{\IDNotation,\IDbNotation}$: the component in $\IDNotation$th row and $\IDbNotation$th column of $\NotationMat \in \mathbb{R}^{\DimANotation \times \DimBNotation}$
	
	\item
	$\VEC{\NotationMat}:=[ 
	\El{\NotationMat}{1,1} , \dots, \El{\NotationMat}{\DimANotation,1} , 
	\El{\NotationMat}{1,2} , \dots, \El{\NotationMat}{\DimANotation,2} , 
	$ $\dots ,
	\El{\NotationMat}{1,\DimBNotation} ,
	\dots, \El{\NotationMat}{\DimANotation,\DimBNotation}	]^{\Mytop} 
	\in \mathbb{R}^{\DimANotation  \DimBNotation} $: 
	the vectorization of $\NotationMat \in \mathbb{R}^{\DimANotation \times \DimBNotation}$

	\item
	$\VECH{\NotationMat}
	:=[ 
	\El{\NotationMat}{1,1} , \dots, \El{\NotationMat}{\DimANotation,1} , 
	\El{\NotationMat}{2,2} , \dots, \El{\NotationMat}{\DimANotation,2} , 
	$ $\dots ,
	\El{\NotationMat}{\DimANotation,\DimANotation} 	]^{\Mytop} 
	\in \mathbb{R}^{\DimANotation ( \DimANotation+1)/2} $: the half vectorization of symmetric $\NotationMat \in \mathbb{R}^{\DimANotation \times \DimANotation}$

	\item
	$\NotationMat_{\mathrm{a}} \otimes \NotationMat_{\mathrm{b}} \in \mathbb{R}^{\DimANotation_{\mathrm{a}} \DimANotation_{\mathrm{b}} \times \DimBNotation_{\mathrm{a}} \DimBNotation_{\mathrm{b}}}$: 
	Kronecker product of $\NotationMat_{\mathrm{a}} \in \mathbb{R}^{\DimANotation_{\mathrm{a}} \times \DimBNotation_{\mathrm{a}}}$ and $\NotationMat_{\mathrm{b}} \in \mathbb{R}^{\DimANotation_{\mathrm{b}} \times \DimBNotation_{\mathrm{b}}}$ defined as
	$$
	\NotationMat_{\mathrm{a}} \otimes \NotationMat_{\mathrm{b}} =
	\begin{bmatrix}
	\El{\NotationMat_{\mathrm{a}}}{1,1} \NotationMat_{\mathrm{b}} & \hdots & \El{\NotationMat_{\mathrm{a}}}{1,\DimBNotation_{\mathrm{a}}} \NotationMat_{\mathrm{b}} \\
	\vdots & \ddots & \vdots \\
	\El{\NotationMat_{\mathrm{a}}}{\DimANotation_{\mathrm{a}},1} \NotationMat_{\mathrm{b}} & \hdots & \El{\NotationMat_{\mathrm{a}}}{\DimANotation_{\mathrm{a}},\DimBNotation_{\mathrm{a}}} \NotationMat_{\mathrm{b}} 
	\end{bmatrix}			
	$$
		
	\item
	$\EigMin{\NotationMat}$: 
	the minimum eigenvalue of symmetric $\NotationMat \in \mathbb{R}^{\DimANotation \times \DimANotation}$
	
	\item 
	${\partial \NotationFunc(\NotationVec)^{\Mytop}}/{\partial \NotationVec} \in \mathbb{R}^{\DimANotation \times \DimBNotation}$: 
	the partial derivative of $\NotationFunc(\NotationVec) \in \mathbb{R}^{\DimBNotation}$ regarding $\NotationVec \in \mathbb{R}^{\DimANotation}$, where 
	$\El{{\partial \NotationFunc(\NotationVec)}/{\partial \NotationVec^{\Mytop}}}{\IDNotation,\IDbNotation}={\partial \El{\NotationFunc(\NotationVec)}{\IDNotation}}/{\partial \El{\NotationVec}{\IDbNotation}}$.	
	
	\item 
	$\Expect{\NotationSto}[\NotationFunc(\NotationSto)]$: 
	the expectation of a function $\NotationFunc(\NotationSto)$ regarding a random variable $\NotationSto$ that follows a probability distribution function (PDF)  denoted by  $\PDF{\StoParam}$.

\end{itemize}

\section{Problem setting}\label{sec_problem}

Sections \ref{sec_problem_target} and \ref{sec_problem_soc} describe the target system and the main problem considered in this study, respectively.

\subsection{Target system}\label{sec_problem_target}

We focus on the following input-affine nonlinear system with an uncertain parameter $\StoParam \in \SetSto$ for a given set $\SetSto \subseteq \mathbb{R}^{\DimSP}$:
\begin{align}
	\dotState(\MyT,\StoParam) 
	&= \Drift\big(\sysx(\MyT,\StoParam),\StoParam\big) 
	+ \InMat\big(\sysx(\MyT,\StoParam),\StoParam\big) \Input(\sysx(\MyT,\StoParam)) 
	,\label{eq:def_sys}
\end{align}
where $\sysx(\MyT,\StoParam) \in \mathbb{R}^{\DimX}$ is the state that depends on $\StoParam$ at the time $\MyT$.
For each $\StoParam \in \SetSto$,
let the drift term $\Drift: \mathbb{R}^{\DimX} \times {\SetSto}\to \mathbb{R}^{\DimX}$,
the input vector field $\InMat: \mathbb{R}^{\DimX} \times {\SetSto} \to  \mathbb{R}^{\DimX \times \DimU}$,
and
a state feedback controller $\Input: \mathbb{R}^{\DimX} \to \mathbb{R}^{\DimU}$ 
be locally Lipschitz continuous in $\sysx(\MyT,\StoParam)$ and satisfy $\Drift(0,\StoParam)=0$.
They obey Assumptions \ref{ass:polypoly_sys} and \ref{ass:StoParam} throughout this paper.

\begin{assumption}[{\MyHighlight{Polytopic nonlinear dynamics}}] \label{ass:polypoly_sys}
	There exist 
	a positive integer $\NumAffine$, 
	continuous functions $\sDrift{\IDAffine}: \SetSto \to  [0,1]$,
	locally Lipschitz continuous functions $\xDrift{\IDAffine}: \mathbb{R}^{\DimX} \to \mathbb{R}^{\DimX}$,
	and
	locally Lipschitz continuous functions $\xInMat{\IDAffine}: \mathbb{R}^{\DimX} \to  \mathbb{R}^{\DimX \times \DimU}$ 
	for $\IDAffine \in \{1,\dots, \NumAffine \}$
	that satisfy the following relations: 	
	\begin{align}
		\Drift(\sysx,\StoParam)
		&=\sum_{\IDAffine=1}^{\NumAffine}  \sDrift{\IDAffine}(\StoParam) \xDrift{\IDAffine}(\sysx)
		,\\
		\InMat(\sysx,\StoParam)
		&=\sum_{\IDAffine=1}^{\NumAffine} \sDrift{\IDAffine}(\StoParam) \xInMat{\IDAffine}(\sysx) 
		,
		\\
		\sum_{\IDAffine=1}^{\NumAffine} \sDrift{\IDAffine}(\StoParam)&=1
		, \label{eq:sum_sDrift_one}
		\\
		\forall \IDAffine,\quad
		\xDrift{\IDAffine}(0)
		&=0
		,	
	\end{align}
	where $\NumAffine$,  $\sDrift{\IDAffine}$, $\xDrift{\IDAffine}$, and $\xInMat{\IDAffine}$  for $\IDAffine \in \{1,\dots, \NumAffine \}$ are known.	
\end{assumption}

While each $\xDrift{\IDAffine}$ is later assumed to be polynomial for robust stability, it is not needed for other technical points.

\begin{assumption}[{\MyHighlight{Stochastic parameter}}] \label{ass:StoParam}
	$\;$	
	\begin{enumerate}
		\item[\MyLabelStoParamUnknown]
		The value of $\StoParam$ is unknown and constant in $\MyT$.  	
		
		\item[\MyLabelPDFknown]
		The PDF $\PDF{\StoParam}$ of $\StoParam$ is known.
	\end{enumerate}
\end{assumption}

\subsection{Stochastic optimal control problem with robust stability}\label{sec_problem_soc}

Our control objective is to minimize the infinite horizon cost function $\costJ(\Input,\IniState,\StoParam)$ under the uncertainty of $\StoParam$:
\begin{align}  
	\costJ(\Input,\IniState,\StoParam)
	&:=
	\lim_{\terminalT \to \infty} 
	\int_{0}^{\terminalT} 
	\Big( 
	\xCost( \sysx(\MyT,\StoParam) ) 
	\nonumber\\&	\quad\quad
	+
	\frac{1}{2} \Input(\sysx(\MyT,\StoParam))^{\Mytop} \uCost( \sysx(\MyT,\StoParam) ) \Input(\sysx(\MyT,\StoParam)) 
	\Big)  
	\mathrm{d}\MyT
	, \label{eq:def_cost}
\end{align}
where $\xCost: \mathbb{R}^{\DimX} \to \mathbb{R}$ is continuous and globally positive definite, that is, $\lim_{\|\sysx\|\to \infty}\xCost(\sysx)\neq 0$.
The function $\uCost: \mathbb{R}^{\DimX} \to \mathbb{R}^{\DimU \times \DimU}$ is continuous and positive definite symmetric for every $\sysx \in \mathbb{R}^{\DimX}$.
The initial state $\sysx(0,\StoParam):=\IniState$ is stochastic and is assumed to obey a PDF $\PDF{\IniState}$ that satisfies the property that the expectation of every monomial of $\IniState$ is finite.

We focus on certain requirements for the design of the feedback controller $\Input$.
The controller $\Input$ cannot utilize the uncertain parameter $\StoParam$ directly because it is unknown under Assumption \ref{ass:StoParam} \MyLabelStoParamUnknown.
Such a controller should aim to minimize  the expectation of $\costJ(\Input,\IniState,\StoParam)$ regarding $\StoParam$.
Offline design of the controller is desirable so that the feedback system can be performed in real time.
Moreover, ensuring stability of the feedback system \eqref{eq:def_sys} is crucial under the uncertainty of $\StoParam  \in \SetSto$.
According to these requirements, this study states the main problem:

{\textbf{Main problem:}} 
Design a feedback controller $\Input$ to minimize the expected cost function and guarantee the robust global asymptotic stability:
	\begin{align}
		&
		\min_{\Input} 
		\Expect{\StoParam}[
		\Expect{\IniState}[    
		\costJ(\Input,\IniState,\StoParam) 
		] 
		]
		,
		\nonumber\\&\quad
		\;\;
		\mathrm{s.t.}
		\;\;
		\forall 
		\StoParam  \in \SetSto
		,\;
		\forall 
		\IniState  \in  \mathbb{R}^{\DimX}
		,\;
		\lim_{\MyT \rightarrow \infty } \sysx(\MyT,\StoParam) = 0 
		.\label{eq:optim_functional}
	\end{align}

\begin{remark}[{\MyHighlight{Difficulties of the main problem}}]\label{rem:difficulties}
The main problem involves the following three difficulties, which are addressed in the next section.
\begin{enumerate}

\item[\MyLabelDifExplicit]
The cost function $\costJ$ and an optimal controller $\Input$ to the problem are not given in explicit form (without the integral form) because of the nonlinearity of the system \eqref{eq:def_sys}.

\item[\MyLabelDifMarkov]
The time-invariant characteristic of $\StoParam$ does not satisfy the Markov property; namely, $\sysx(\MyT_{2},\StoParam)$ depends not only on the value of $\sysx(\MyT_{1},\StoParam)$ but also on $\StoParam$ for $\MyT_{1} < \MyT_{2}$.
Efficient tools for optimal control, such as HJB equations and the principle of optimality (dynamic programming), are not applicable because they are based on the Markov property.

\item[\MyLabelDifStability]
The controller needs to realize not only the optimality but also the robust global asymptotic stability of the nonlinear system \eqref{eq:def_sys}.
Ensuring the stability is challenging because of both the nonlinearity and uncertainty.
	
\end{enumerate}

\end{remark}

\section{Proposed method: Stochastic optimal control with robust stability}\label{sec_method}

This section presents the proposed method to solve the main problem.
Section \ref{sec_overview} provides an overview of the method.
The technical points for establishing the method are detailed in Sections \ref{sec_parametric_approx}--\ref{sec_implementation}.

\subsection{Overview of the proposed method}\label{sec_overview}

To overcome the three difficulties {\MyLabelDifExplicit}--{\MyLabelDifStability} in Remark \ref{rem:difficulties} associated with the main problem, 
we establish the proposed method based on three techniques: parametric approximation, a gradient method, and SOS. 
The main problem \eqref{eq:optim_functional} is relaxed into the following minimization problem:
\begin{align}
&
\min_{\InVecParam,\LyapMat,\ASparam} 
\ObjF(\InVecParam,\LyapMat,\ASparam)
,\label{eq:optim_param}
\\
\ObjF(\InVecParam,\LyapMat,\ASparam)
&:=
\Expect{\StoParam}[\Expect{\IniState}[    
\EachValue(\IniState,\VParam(\InVecParam,\StoParam))
] ]
+
\PenaltyFuncRegu(\InVecParam, \LyapMat,\ASparam)
,\label{eq:def_ObjF}
\end{align}
where $\EachValue$ denotes an approximate of the true cost $\costJ$, which is characterized by a parameter function $\VParam(\InVecParam,\StoParam)$.
The symbol $\InVecParam$ is a decision variable that characterizes the controller $\Input$.
The other decision variables $(\LyapMat,\ASparam)$ and the penalty function $\PenaltyFuncRegu$ are later introduced to guarantee the stability.
An overview of the proposed method using the three techniques is described below.

Firstly, the parametric approximation is utilized to express the controller and cost function in an explicit form, which addresses the difficulty {\MyLabelDifExplicit}.
The controller $\Input(\sysx)$ and the cost function $\costJ(\Input,\sysx,\StoParam)$ in \eqref{eq:optim_functional} are replaced by parametric functions $\InpOpt(\sysx, \InVecParam) $ and $\EachValue(\sysx, \VParam(\InVecParam,\StoParam))$ in \eqref{eq:optim_param}, respectively, as follows:
\begin{align}
\Input(\sysx)
&=
\InpOpt(\sysx, \InVecParam) 
,\label{eq:apply_FBcontrol}
\\
\costJ(\InpOpt(\NonArg, \InVecParam), \sysx,\StoParam)
&=
\EachValue(\sysx, \VParam(\InVecParam,\StoParam))
+ \ErrV(\sysx,\StoParam,\InVecParam) 
,
\end{align}
where $\ErrV(\sysx,\StoParam,\InVecParam) $ denotes the approximation residual.
The functions $\EachValue(\sysx, \VParam(\InVecParam,\StoParam))$ and $\InpOpt(\sysx, \InVecParam)$ are parameterized by  $\VParam(\InVecParam,\StoParam) \in \mathbb{R}^{\DimVParam}$ and $\InVecParam \in \mathbb{R}^{\DimInVecParam}$, respectively, with basis functions of $\sysx$.
It should be noted that $\VParam(\InVecParam,\StoParam)$ depends on $\InVecParam$ and $\StoParam$ because $\costJ(\InpOpt(\NonArg, \InVecParam), \sysx,\StoParam)$ depends on them.
The details of this parametric approximation are presented in Section \ref{sec_parametric_approx}.

Secondly, a gradient-based approach is efficient for overcoming the difficulty {\MyLabelDifMarkov} because the Markov property is not required, unlike the principle of optimality and HJB equations.
The gradient method is used to obtain a suboptimal solution to \eqref{eq:optim_param}:
\begin{align}
{\Ite{\InVecParam}{\IDIteGrad+1}}
&= 
{\Ite{\InVecParam}{\IDIteGrad}}
- {\CoefGrad{\IDIteGrad}} \frac{\partial }{\partial \InVecParam } 
\ObjF(  {\Ite{\InVecParam}{\IDIteGrad}}   ,{\Ite{\LyapMat}{\IDIteGrad}} ,{\Ite{\ASparam}{\IDIteGrad}})
, \label{eq:gradient_updatate_InVecParam}
\\
{\Ite{\LyapMat}{\IDIteGrad+1}}
&= 
{\Ite{\LyapMat}{\IDIteGrad}}
- {\CoefGrad{\IDIteGrad}} \frac{\partial }{\partial \LyapMat } 
\PenaltyFuncRegu(  {\Ite{\InVecParam}{\IDIteGrad}}   ,{\Ite{\LyapMat}{\IDIteGrad}} ,{\Ite{\ASparam}{\IDIteGrad}})
, \label{eq:gradient_updatate_LyapMat}
\\
{\Ite{\ASparam}{\IDIteGrad+1}}
&= 
{\Ite{\ASparam}{\IDIteGrad}}
- {\CoefGrad{\IDIteGrad}} \frac{\partial }{\partial \ASparam } 
\PenaltyFuncRegu(  {\Ite{\InVecParam}{\IDIteGrad}}   ,{\Ite{\LyapMat}{\IDIteGrad}} ,{\Ite{\ASparam}{\IDIteGrad}})
, \label{eq:gradient_updatate_ASparam}
\end{align}
where the superscript ${\Ite{}{\IDIteGrad}}$ denotes the decision variables at the $\IDIteGrad$th iteration.
After iterations, a suboptimal solution and the corresponding feedback controller $\InpOpt(\sysx, {\Ite{\InVecParam}{\IDIteGrad}})$ are obtained.
The details of the gradient derivation are presented in Section \ref{sec_gradient}.

Thirdly, the penalty function $\PenaltyFuncRegu(\InVecParam,\LyapMat,\ASparam)$ in \eqref{eq:def_ObjF} helps the controller ensure the robust global asymptotic stability, which is related to the difficulty {\MyLabelDifStability}.
The penalty function is formulated based on SOS that is a promising tool for guaranteeing stability of polynomial nonlinear systems.
The details of the penalty function are presented in Section \ref{sec_SOS}.

Sections \ref{sec_parametric_approx}--\ref{sec_SOS} describe the technical aspects of the proposed method.
Finally, Section \ref{sec_implementation} describes the implementation of the proposed method.

\begin{remark}
		The proposed gradient-based approach is an advanced version of our previous method only for linear systems \citep{Ito24Cyb,ItoCDC16}.
		Although the previous method has not been applied to nonlinear systems due to the difficulties {\MyLabelDifExplicit} and {\MyLabelDifStability}, this study overcomes them, as described above.
\end{remark}

\subsection{Parametric approximation} \label{sec_parametric_approx}

This subsection presents the details of the parameterization of the controller $\InpOpt(\sysx, \InVecParam)$ and the approximate cost function $\EachValue(\sysx,\VParam(\InVecParam,\StoParam))$.
These are defined as follows:
\begin{align}
\InpOpt(\sysx, \InVecParam) 
&:= 
\uBasis(\sysx)^{\Mytop} \InVecParam
,  \label{eq:def_parametric_input}
\\ 
\EachValue(\sysx,\VParam(\InVecParam,\StoParam)) 
&:=
\JBasis(\sysx)^{\Mytop}  \VParam(\InVecParam,\StoParam)
, \label{eq:def_parametric_cost}
\end{align}
where $\uBasis: \mathbb{R}^{\DimX} \to \mathbb{R}^{\DimInVecParam \times \DimU}$ and $\JBasis: \mathbb{R}^{\DimX} \to \mathbb{R}^{\DimVParam}$ are a locally Lipschitz continuous basis function and a $C^{1}$ continuous basis function satisfying $\JBasis(0)=0$, respectively.
The controller parameter $\InVecParam \in \mathbb{R}^{\DimInVecParam}$ is optimized as shown in \eqref{eq:gradient_updatate_InVecParam}. 
For each $\InVecParam$ and $\StoParam$, the parameter $\VParam(\InVecParam,\StoParam) \in \mathbb{R}^{\DimVParam}$ should be determined such that 
the function $\EachValue(\sysx,\VParam(\InVecParam,\StoParam))$ in \eqref{eq:def_parametric_cost} approximates
the true cost function $\costJ(\InpOpt(\NonArg, \InVecParam), \sysx,\StoParam)$, namely,
a distance between the two functions should be minimized.
Introducing such a distance is challenging because $\costJ(\InpOpt(\NonArg, \InVecParam), \sysx,\StoParam)$ cannot be obtained directly.

To overcome this challenge, 
we define an optimal parameter  $\VParam(\InVecParam,\StoParam)$ as a minimizer to the following distance based on the Bellman residual $\Bellman(\sysx,\preVParam,\InVecParam,\StoParam)$:
\begin{align}
	\VParam(\InVecParam,\StoParam)
	&\in
	\argmin_{\preVParam \in \mathbb{R}^{\DimVParam}}	
	\Big(
	\MarginalMap[\Big]{   \Bellman( \NonArg, \preVParam,\InVecParam, \StoParam)^{2}  }
	+
	{\CoefReguVParam} \|\preVParam\|^{2} 
	\Big)	
	,\label{eq:def_OptVParam}
	\\
		\Bellman(\sysx,\preVParam,\InVecParam,\StoParam)
		&:=
		\xCost(\sysx)
		+
		\frac{1}{2}
		\InpOpt(\sysx, \InVecParam)^{\Mytop}
		\uCost(\sysx)
		\InpOpt(\sysx, \InVecParam)
		\nonumber\\&\quad
		+
		\frac{\partial \EachValue(\sysx,\preVParam) }{\partial \sysx}^{\Mytop}
		\Big(
		\Drift(\sysx,\StoParam) 
		+
		\InMat(\sysx,\StoParam)
		\InpOpt(\sysx, \InVecParam)
		\Big)
		,\label{eq:Bellman_eq}
\end{align}
where  $\CoefReguVParam\geq 0$ denotes a free parameter.
The linear operator $\MarginalMap[]{\NonArg }$ is given by
\begin{align}
		\MarginalMap[\Big]{ \NotationFunc(\NonArg) } 
		&:=
		\int_{\mathbb{R}^{\DimX}} \MarginalWeight(\sysx) \NotationFunc(\sysx) \mathrm{d} \sysx
		,\label{eq:def_MarginalMap}
\end{align}
for any continuous function $ \NotationFunc:{\mathbb{R}^{\DimX}} \to \mathbb{R}^{\DimANotation\times \DimBNotation} $ with any $(\DimANotation,\DimBNotation)$, where $\MarginalWeight: {\mathbb{R}^{\DimX}} \to [ 0, \infty] $ is a predefined weight function such that $\MarginalMap[]{\NonArg }$ is linear, bounded, and well-defined for every $ \NotationFunc$.
In the following, we present Theorem \ref{thm:BellmanResidual} along with Remark \ref{rem:BellmanResidual} to justify that  \eqref{eq:def_OptVParam} is  appropriate  for approximating $\EachValue(\sysx,\VParam(\InVecParam,\StoParam))$.
Next, we derive Theorem \ref{thm:explicit_cost} to obtain a solution  $\VParam(\InVecParam,\StoParam)$ to \eqref{eq:def_OptVParam} in explicit form.

\begin{theorem}[{\MyHighlight{Bellman residual}}]\label{thm:BellmanResidual}
Given parameters $\preVParam \in \DimVParam$, $\InVecParam \in \DimInVecParam$, and $\StoParam \in \SetSto$, 
suppose 
that the feedback system \eqref{eq:def_sys} with $\Input(\sysx)=\InpOpt(\sysx, \InVecParam)$ is globally asymptotically stable\footnote{%
For the given $\InVecParam \in \DimInVecParam$, and $\StoParam \in \SetSto$, the global asymptotic stability implies that for every $\IniState  \in  \mathbb{R}^{\DimX}$, we have 
$\sup_{\MyT \geq 0 } \|\sysx(\MyT,\StoParam)\| < \infty$ and $\lim_{\MyT \rightarrow \infty } \sysx(\MyT,\StoParam) = 0$.}
and 
that for each $\sysx \in \mathbb{R}^{\DimX}$, $\costJ(\InpOpt(\NonArg, \InVecParam),\sysx,\StoParam) < \infty $ holds.
If there exists $\BellResiCoef\geq 0$ such that
\begin{align}
&\forall \sysx \in  \mathbb{R}^{\DimX},
\nonumber\\
&\Bellman(\sysx,\preVParam,\InVecParam,\StoParam)^{2}
\leq
\BellResiCoef^{2}
\Big(
\xCost(\sysx)
+
\frac{1}{2}
\InpOpt(\sysx, \InVecParam)^{\Mytop}
\uCost(\sysx)
\InpOpt(\sysx, \InVecParam)
\Big)^{2}
,\label{eq:Bellman_residual_bound}
\end{align}	
then 
we have 
\begin{align}
\forall \sysx \in \mathbb{R}^{\DimX}
,
|
\EachValue(\sysx,\preVParam)
-
\costJ(\InpOpt(\NonArg, \InVecParam), \sysx,\StoParam)
|
\leq
\BellResiCoef
\costJ(\InpOpt(\NonArg, \InVecParam), \sysx,\StoParam) 
.\label{eq:Bellman_eq_equivalent_err} 		
\end{align}	
\end{theorem}
\begin{\MyProof}
The proof is described in Appendix \ref{pf:BellmanResidual}.
\end{\MyProof}
\begin{remark}[{\MyHighlight{Contributions of  Theorem \ref{thm:BellmanResidual}}}]\label{rem:BellmanResidual}
Theorem \ref{thm:BellmanResidual} indicates that 
$\BellResiCoef$ characterizes the relation between the Bellman residual $\Bellman(\sysx,\preVParam,\InVecParam,\StoParam)$
and the approximation error 
$
|
\EachValue(\sysx,\preVParam)
-
\costJ(\InpOpt(\NonArg, \InVecParam), \sysx,\StoParam)
|
$.
Specifically, if $\Bellman(\sysx,\preVParam,\InVecParam,\StoParam)$ is bounded with a small $\BellResiCoef$,
the approximation error decreases.
This relation justifies minimizing the norm of the Bellman residual $\MarginalMap{   \Bellman( \NonArg, \preVParam,\InVecParam, \StoParam)^{2}  }$ that constitutes the first term in  \eqref{eq:def_OptVParam}.
The second term in \eqref{eq:def_OptVParam} ensures that a solution  $\VParam(\InVecParam,\StoParam)$ to \eqref{eq:def_OptVParam} is uniquely determined.
\end{remark}

\begin{theorem}[{\MyHighlight{Solution to the optimal approximation}}]\label{thm:explicit_cost}
Suppose that Assumption \ref{ass:polypoly_sys} holds.
Suppose that $\CoefReguVParam>0$,
$
\MarginalMap[]{ 
 \| {\xBterm} \|^{2}
}
< \infty
$,
and
$
\MarginalMap[]{ 
\| {\xAterm{\IDAffine}}\|^{2} 
}
< \infty
$ for $\IDAffine \in \{1,\dots, \NumAffine \}$
hold, where
\begin{align}
\xBterm
&:=
\begin{bmatrix}
	\xCost(\sysx) 
	\\ 
	{\VEC[\big]{  
			\uBasis(\sysx)	\uCost(\sysx)	\uBasis(\sysx)^{\Mytop}/2
	}}
\end{bmatrix}
,\\
{\xAterm{\IDAffine}}
&:=
{\VEC[\Big]{  
		\frac{\partial \JBasis(\sysx)}{\partial \sysx^{\Mytop}}
		[ \xDrift{\IDAffine}(\sysx)  ,  \xInMat{\IDAffine}(\sysx) \uBasis(\sysx)^{\Mytop}   ]
}}
. 	
\end{align}
For any $\InVecParam  \in \mathbb{R}^{\DimInVecParam}$ and any $\StoParam \in \SetSto$, 
the solution $\VParam(\InVecParam,\StoParam)$ to \eqref{eq:def_OptVParam} is uniquely determined by
\begin{align}
	&
	\VParam(\InVecParam,\StoParam)
	\nonumber\\&
	=
	-
	\Big(
	\sum_{\IDAffine=1}^{\NumAffine}
	\sum_{\IDbAffine=1}^{\NumAffine}
	\sDrift{\IDAffine}(\StoParam) 
	\sDrift{\IDbAffine}(\StoParam)
	{\newDDmat{\InVecParam}{\IDAffine}{\IDbAffine}} 
	\Big)^{-1}  
	\sum_{\IDAffine=1}^{\NumAffine}
	\sDrift{\IDAffine}(\StoParam)
	{\newDNmat{\InVecParam}{\IDAffine}} 	
	, \label{eq:solution_OptVParam}
\end{align}		
where
\begin{align}
	&{\newDDmat{\InVecParam}{\IDAffine}{\IDbAffine}} 
	\nonumber\\&
	:=
	\Big( 
	\begin{bmatrix}
		1 \\ \InVecParam
	\end{bmatrix}
	\otimes  \Identity{\DimVParam} 
	\Big)^{\Mytop} 
	\MarginalMap[\Big]{ 
		{\xAterm{\IDAffine}}
		{\xAterm{\IDbAffine}}^{\Mytop}
	}
	\Big( 
	\begin{bmatrix}
		1 \\ \InVecParam
	\end{bmatrix}
	\otimes  \Identity{\DimVParam} 
	\Big)
	+
	{\CoefReguVParam}
	\Identity{\DimVParam}
	,\\
	&{\newDNmat{\InVecParam}{\IDAffine}} 
	:=
	\Big( 
	\begin{bmatrix}
		1 \\ \InVecParam
	\end{bmatrix}
	\otimes  \Identity{\DimVParam} 
	\Big)^{\Mytop}  
	\MarginalMap[\Big]{ 
		{\xAterm{\IDAffine}}
		{\xBterm}^{\Mytop}
	}
	\begin{bmatrix}
		1 \\ \InVecParam \otimes \InVecParam
	\end{bmatrix}
.	
\end{align}
Furthermore, the assumption of $\CoefReguVParam>0$ is not needed if ${\newDDmat{\InVecParam}{\IDAffine}{\IDbAffine}} \succ 0$ holds.
\end{theorem}
\begin{\MyProof}
The proof is described in Appendix \ref{pf:explicit_cost}.
\end{\MyProof}
\begin{remark}[{\MyHighlight{Contributions of  Theorem \ref{thm:explicit_cost}}}]
	Theorem \ref{thm:explicit_cost} gives an explicit solution  $\VParam(\InVecParam,\StoParam)$ to \eqref{eq:def_OptVParam}.
	Substituting this solution into  \eqref{eq:def_parametric_cost} yields an explicit form of the approximate cost function $\EachValue(\IniState,\VParam(\InVecParam,\StoParam))=
	\JBasis(\sysx)^{\Mytop}  \VParam(\InVecParam,\StoParam)
	$.
	This cost is an optimal approximate of the true cost $\costJ(\InpOpt(\NonArg, \InVecParam), \sysx,\StoParam)$ in the sense that the Bellman-based distance in \eqref{eq:def_OptVParam} is minimized.
\end{remark}

\subsection{The gradient of the approximated cost function}  \label{sec_gradient}

The gradient of the approximate cost function $\EachValue(\IniState,\VParam(\InVecParam,\StoParam))$ is required 
to realize the proposed gradient method in \eqref{eq:gradient_updatate_InVecParam}.
We derive the gradient below.

\begin{theorem}[{\MyHighlight{Explicit gradient of the  cost}}]\label{thm:explicit_grad}
Suppose that all the conditions assumed in Theorem \ref{thm:explicit_cost} hold. 
For any $\InVecParam \in \mathbb{R}^{\DimInVecParam}$ and any $\StoParam \in \SetSto$, the gradient of the approximate cost $\EachValue(\IniState,\VParam(\InVecParam,\StoParam))$  with \eqref{eq:solution_OptVParam} is given by
	\begin{align}
		&
		\frac{\partial }{\partial \El\InVecParam{\IDbEl}}
		\Expect{\IniState}[    
		\EachValue(\IniState,\VParam(\InVecParam,\StoParam))
		]  
		\nonumber\\&
		=
		-
		\Expect{\IniState}[   
		\JBasis(\IniState) 
		]^{\Mytop} 
		\Big(
		\sum_{\IDAffine=1}^{\NumAffine}
		\sum_{\IDbAffine=1}^{\NumAffine}
		\sDrift{\IDAffine}(\StoParam) 
		\sDrift{\IDbAffine}(\StoParam)
		{\newDDmat{\InVecParam}{\IDAffine}{\IDbAffine}} 
		\Big)^{-1} 
		\nonumber
		\\& \quad \times
		\sum_{\IDAffine=1}^{\NumAffine}
		\sDrift{\IDAffine}(\StoParam)
		\Big(
		\frac{\partial    {\newDNmat{\InVecParam}{\IDAffine}}    }{\partial \El\InVecParam{\IDbEl}}
		+
		\sum_{\IDbAffine=1}^{\NumAffine}
		\sDrift{\IDbAffine}(\StoParam)
		\frac{\partial {\newDDmat{\InVecParam}{\IDAffine}{\IDbAffine}}  }{\partial \El\InVecParam{\IDbEl}}
		\VParam(\InVecParam,\StoParam))
		\Big)
		.\label{eq:explicit_cost_derivative}	
	\end{align}	
Furthermore, the assumption of $\CoefReguVParam>0$ is not needed if ${\newDDmat{\InVecParam}{\IDAffine}{\IDbAffine}} \succ 0$  holds in a neighborhood of $\InVecParam$.
\end{theorem}
\begin{\MyProof}
The proof is described in Appendix \ref{pf:explicit_grad}.
\end{\MyProof}
\begin{remark}[{\MyHighlight{Contributions of  Theorem \ref{thm:explicit_grad}}}]
Theorem \ref{thm:explicit_grad} provides the gradient of the approximate cost function $\EachValue(\IniState,\VParam(\InVecParam,\StoParam))$ in explicit form.
This gradient enables implementation of the gradient method in \eqref{eq:gradient_updatate_InVecParam} to solve the minimization problem \eqref{eq:optim_param}. 
\end{remark}

\subsection{SOS-based stabilization combined with the gradient method} \label{sec_SOS}

In this subsection, the robust global asymptotic stability of the feedback system \eqref{eq:def_sys} with $\Input(\sysx)=\InpOpt(\sysx, \InVecParam)$ is ensured while the gradient method optimizes the controller.
Our strategy for guaranteeing the stability is based on SOS.
For a monomial basis function $ \SOSBasis(\sysx)\in \mathbb{R}^{\DimSOSBasis}$, a Lyapunov function $\LyapF(\sysx) $ and its time derivative $\dot{\LyapF}(\sysx,\StoParam)$ are defined by the following SOS forms:
\begin{align}
	\LyapF(\sysx) 
	&:=
	\SOSBasis(\sysx)^{\Mytop} \LyapMat \SOSBasis(\sysx)
	,\label{eq:def_Lyap}
	\\
	\dot{\LyapF}(\sysx(\MyT,\StoParam),\StoParam) 
	&= 	- 2 \sum_{\IDAffine=1}^{\NumAffine}  \sDrift{\IDAffine}(\StoParam) 
	\SOSBasis(\sysx(\MyT,\StoParam))^{\Mytop} 
	\nonumber\\&
	\quad\quad\quad\times
	{\SOSdotVpreMat{\IDAffine}}(\sysx(\MyT,\StoParam),\LyapMat,\InVecParam)
	\SOSBasis(\sysx(\MyT,\StoParam))
	,\label{eq:def_dLyap_dt}
\end{align}
where  $\LyapMat \in  \mathbb{R}^{\DimSOSBasis \times \DimSOSBasis}$ and  ${\SOSdotVpreMat{\IDAffine}}(\sysx,\LyapMat,\InVecParam) \in  \mathbb{R}^{\DimSOSBasis \times \DimSOSBasis}$ are the decision variable in the minimization \eqref{eq:optim_param} and a function to be clarified later, respectively.
The robust global asymptotic stability holds if for every $\StoParam$, $\LyapF(\sysx) $ and $-\dot{\LyapF}(\sysx,\StoParam)$ are globally positive definite and  $\LyapF(\sysx) $ is radially unbounded, that is, $\LyapF(\sysx)>0$ and $\dot{\LyapF}(\sysx,\StoParam)<0$ for $\sysx \neq 0$,  $\LyapF(0)=\dot{\LyapF}(0,\StoParam)=0$,  $\lim_{\|\sysx\| \to \infty} \LyapF(\sysx) = \infty $, and $\lim_{\|\sysx\| \to \infty} \dot{\LyapF}(\sysx,\StoParam) \neq 0 $.
While this strategy is based on groundbreaking results \citep[Theorem 1]{Xu09TAC} and \citep[Theorem 6]{Prajna04}, we have novelty and advantages compared with them to be discussed in Remark \ref{rem:related_work} at the end of this subsection.

The goal of this subsection is to design a penalty function $\PenaltyFuncRegu(\InVecParam,\LyapMat,\ASparam)$ and initial decision variables  $(  {\Ite{\InVecParam}{0}}   ,{\Ite{\LyapMat}{0}} ,{\Ite{\ASparam}{0}})$ so that the aforementioned SOS-based strategy is realized in the gradient method in \eqref{eq:gradient_updatate_InVecParam}--\eqref{eq:gradient_updatate_ASparam}. 
We present key theorems to justify the design in the following. 
Based on Lemma \ref{thm:design_penaF}, we introduce the details of 
the penalty function $\PenaltyFuncRegu(\InVecParam,\LyapMat,\ASparam)$ and the decision variable $\ASparam$. 
Theorem \ref{thm:stabilizing_optim} implies that the proposed penalty function guarantees the robust global asymptotic stability.
Theorem \ref{thm:stabilizing_initial} provides appropriate initial values $(  {\Ite{\InVecParam}{0}}   ,{\Ite{\LyapMat}{0}} ,{\Ite{\ASparam}{0}})$ used at the beginning of the gradient method. 
The following definitions and assumptions of monomial bases are used.

\begin{definition}[{\MyHighlight{Monomial bases}}]\label{def:basis}
	A function $\DEFSOSBasis: \mathbb{R}^{\DimX} \to \mathbb{R}^{\DEFDimASOSBasis \times \DEFDimBSOSBasis}$ is said to be a monomial basis 
	if each component of $\DEFSOSBasis(\sysx)$ is a primitive monomial $\prod_{\IDEl=1}^{\DimX} \El{\sysx}{\IDEl}^{\El{\degPow}{\IDEl}}$ for some $\degPow \in \{0,1,2,\dots\}^{\DimX}$.
	A monomial basis $\DEFSOSBasis(\sysx)$ is said to be strict 
	if $\DEFSOSBasis(\sysx)=0 \Leftrightarrow \sysx=0$ is satisfied.
	A monomial basis $\DEFSOSBBasis(\sysx)$ is said to be the non-redundant form of a monomial basis $\DEFSOSBasis(\sysx)$ if
	all the components of $\DEFSOSBBasis(\sysx)$ consist of those of $\DEFSOSBasis(\sysx)$ and
	are distinct monomials, that is, for every $\IDEl\neq \IDbEl$, ${\El{\DEFSOSBBasis(\sysx)}{\IDEl}}\neq {\El{\DEFSOSBBasis(\sysx)}{\IDbEl}}$ holds for some $\sysx$. 	
\end{definition}

\begin{remark}[{\MyHighlight{Introducing basis functions}}]\label{rem:monomial_bases}
Let
$\SOSBasis: \mathbb{R}^{\DimX} \to \mathbb{R}^{\DimSOSBasis}$ and
$\SOSmatBasis: \mathbb{R}^{\DimX} \to \mathbb{R}^{\DimSOSmatBasis}$ be monomial bases such that $\SOSBasis(\sysx)$ is strict and ${\El{\SOSmatBasis(\sysx)}{1}}=1$.
Let $\SOSdotVBasis: \mathbb{R}^{\DimX} \to \mathbb{R}^{\DimSOSdotVBasis}$ be the non-redundant form of $\SOSmatBasis(\sysx)\otimes \SOSBasis(\sysx) \in \mathbb{R}^{\DimSOSBasis \DimSOSmatBasis }$.
Let $\InBasis:  \mathbb{R}^{\DimX} \to \mathbb{R}^{\DimU \times \DimRowInBasis  }$ be a continuous function such that 
$
{\PolyInMat{\IDAffine}}(\sysx) 
\InBasis(\sysx)
$ is polynomial.
Further guidelines for designing these bases correspond to the following assumptions.
\end{remark}

\begin{assumption}[{\MyHighlight{SOS settings}}] \label{ass:SOS}
The bases $\SOSBasis(\sysx)$, $\SOSmatBasis(\sysx)$, and $\InBasis(\sysx)$ are designed such that the following conditions hold:	
\begin{enumerate}
	\item[\MyLabelFpoly] 
	For each $\IDAffine \in \{1,2,\dots,\NumAffine \}$,
	there exist known polynomial matrices ${\PolyDriftMat{\IDAffine}}: \mathbb{R}^{\DimX} \to \mathbb{R}^{\DimX \times \DimSOSBasis}$  satisfying
	\begin{align}
		\xDrift{\IDAffine}(\sysx) = {\PolyDriftMat{\IDAffine}}(\sysx) \SOSBasis(\sysx)
		. \label{eq:def_PolyDriftMat}
	\end{align}	
	
	\item[\MyLabelUpoly]
	The basis function $\uBasis(\sysx)$ used in the controller $\InpOpt(\sysx, \InVecParam) $ in \eqref{eq:def_parametric_input} is set to
	\begin{align}
		\uBasis(\sysx)&= \SOSBasis(\sysx) \otimes \InBasis(\sysx)^{\Mytop} 
		. \label{eq:def_InBasis}
	\end{align}
	
	\item[\MyLabeldotVpoly]	
	For any  $\LyapMat \in \mathbb{R}^{\DimSOSBasis\times \DimSOSBasis}$, any controller parameter $\InVecParam \in   \mathbb{R}^{\DimInVecParam} $,
	and any $\IDAffine\in\{1,2,\dots,\NumAffine\}$,
	the following function:
	\begin{align}
		&{\SOSdotVpreMat{\IDAffine}}(\sysx,\LyapMat,\InVecParam)
		\nonumber\\&
		:=
		-\frac{  \LyapMat + \LyapMat^{\Mytop}  }{2}
		\frac{\partial \SOSBasis(\sysx)}{\partial \sysx^{\Mytop} }
		(
		{\PolyDriftMat{\IDAffine}}(\sysx) 
		+
		{\PolyInMat{\IDAffine}}(\sysx) 
		\InBasis(\sysx)
		{\InvVEC{\InVecParam}{\DimRowInBasis}{\DimSOSBasis}}
		)		
		\label{eq:def_SOSdotVpreMat}
	\end{align}
	is a polynomial matrix consisting of the basis $\SOSmatBasis(\sysx)\otimes \SOSmatBasis(\sysx)$,
	that is, 
	each $(\IDEl,\IDbEl)$th component satisfies  ${\El{ {\SOSdotVpreMat{\IDAffine}}(\sysx,\LyapMat,\InVecParam) }{\IDEl,\IDbEl}}={\coeffdotVpreSqrtMat{\IDEl}{\IDbEl}}^{\Mytop}  ( \SOSmatBasis(\sysx)\otimes \SOSmatBasis(\sysx) ) $ for some ${\coeffdotVpreSqrtMat{\IDEl}{\IDbEl}} \in \mathbb{R}^{\DimSOSmatBasis^{2}}$,
	where
	${\InvVEC{\InVecParam}{\DimRowInBasis}{\DimSOSBasis}}\in \mathbb{R}^{\DimRowInBasis \times \DimSOSBasis}$ denotes the inverse map of $\VEC{\NonArg}$ satisfying 
	$\InVecParam=\VEC{ {\InvVEC{\InVecParam}{\DimRowInBasis}{\DimSOSBasis}} }$.

	\item[\MyLabelSOS]
	There exist a positive definite symmetric  $\LyapMat \succ 0$ and $\InVecParam \in   \mathbb{R}^{\DimInVecParam} $
	such that
	for every $\IDAffine\in\{1,2,\dots,\NumAffine\}$,
	${\SOSdotVpreMat{\IDAffine}}(\sysx,\LyapMat,\InVecParam)
	+
	{\SOSdotVpreMat{\IDAffine}}(\sysx,\LyapMat,\InVecParam)^{\Mytop} 
	$ is a strict SOS\footnote{A strict SOS is a strict version of SOS \citep[Sec. 2.1]{Scherer06}. The SOS and strict SOS imply ${\SOSpdMat{\IDAffine}} \succeq 0$ and ${\SOSpdMat{\IDAffine}} \succ 0$ in \eqref{eq:SOS_ass}, respectively.} regarding $\SOSmatBasis(\sysx)$,	
	that is, there exists a positive definite symmetric matrix ${\SOSpdMat{\IDAffine}}\succ 0$ satisfying
	\begin{align}
	&{\SOSdotVpreMat{\IDAffine}}(\sysx,\LyapMat,\InVecParam)
	+
	{\SOSdotVpreMat{\IDAffine}}(\sysx,\LyapMat,\InVecParam)^{\Mytop} 
	\nonumber\\&
	=
	(\SOSmatBasis(\sysx) \otimes {\Identity{\DimSOSBasis}} )^{\Mytop}
	{\SOSpdMat{\IDAffine}} 
	(\SOSmatBasis(\sysx) \otimes {\Identity{\DimSOSBasis}} )
	.\label{eq:SOS_ass}
	\end{align}

\end{enumerate}
\end{assumption}

Based on these assumptions, 
we propose the penalty function $\PenaltyFuncRegu:
\mathbb{R}^{\DimInVecParam} \times \mathbb{R}^{\DimSOSBasis\times \DimSOSBasis} \times  \mathbb{R}^{\DimASparam \times \NumAffine } \to \mathbb{R}$ as the following logarithmic barrier function:
\begin{align}
	&
	\PenaltyFuncRegu(\InVecParam,\LyapMat,\ASparam)
	\nonumber\\&
	:=
	\begin{cases}
		-
		\penaW\ln \Big( \det((\LyapMat+\LyapMat^{\Mytop})/2)	
		\prod_{\IDAffine=1}^{\NumAffine}
		\det({\SOSdotVMat{\IDAffine}}(\InVecParam, \LyapMat,\ASparam) )	
		\Big)
		\\
		\qquad\qquad\qquad\qquad
		(\LyapMat+\LyapMat^{\Mytop} \succ 0, \;\;{\SOSdotVMat{\IDAffine}}(\InVecParam, \LyapMat,\ASparam)  \succ 0),
		\\
		\iniPenaltyFuncRegu
		\qquad\qquad\qquad\;\;
		(\mathrm{otherwise}),
	\end{cases}
	\label{eq:def_PenaltyFuncRegu}
\end{align}
where $\penaW>0$ and $\iniPenaltyFuncRegu\in \mathbb{R}$ are constants, and the functions
${\SOSdotVMat{\IDAffine}}(\InVecParam, \LyapMat,\ASparam)$ for $\IDAffine\in \{1,2,\dots,\NumAffine\}$ are defined as follows.
Letting $\ASparam:=[{\ASeach{1}},{\ASeach{2}},\dots,{\ASeach{\NumAffine}}]$ with ${\ASeach{\IDAffine}} \in  \mathbb{R}^{\DimASparam } $, we define ${\SOSdotVMat{\IDAffine}}(\InVecParam, \LyapMat,\ASparam)$ as a unique solution of ${\SOSdotVMat{\IDAffine}}$ satisfying the following conditions:
\begin{align}
\SOSdotVBasis(\sysx)^{\Mytop}
{\SOSdotVMat{\IDAffine}}
\SOSdotVBasis(\sysx)
&=
\SOSBasis(\sysx)^{\Mytop}
{\SOSdotVpreMat{\IDAffine}}(\sysx,\LyapMat,\InVecParam)
 \SOSBasis(\sysx)
,\;
\forall \sysx \in \mathbb{R}^{\DimX}
,\label{eq:def_dotLyapF_new}	
\\
\AScoefmat
{\VECH[\big]{  
		{\SOSdotVMat{\IDAffine}}
}}
&=\AScoefscalar {\ASeach{\IDAffine}}
,\label{eq:def_ASconstraint}	
\\
{\SOSdotVMat{\IDAffine}}&={\SOSdotVMat{\IDAffine}}^{\Mytop}
,\label{eq:def_symmetricSOSdotVMat}	
\end{align}
for some matrix $\AScoefmat$ and some scalar $\AScoefscalar$.
The uniqueness  of ${\SOSdotVMat{\IDAffine}}(\InVecParam, \LyapMat,\ASparam)$ is ensured by the following lemma.

\begin{lemma}[{\MyHighlight{Uniqueness of ${\SOSdotVMat{\IDAffine}}(\InVecParam, \LyapMat,\ASparam)$}}]\label{thm:design_penaF}
Suppose that Assumptions \ref{ass:polypoly_sys} and \ref{ass:SOS} hold.
There exist  $\AScoefmat \in \mathbb{R}^{\DimASparam \times \DimSOSdotVBasis(\DimSOSdotVBasis+1)/2 } $ (for some $\DimASparam\in \{1,2,\dots\}$) and $\AScoefscalar\in \{0,1\}$ that satisfy the following two statements.
For any $\IDAffine$ and any $(\InVecParam,\LyapMat,\ASparam) 
	\in  \mathbb{R}^{\DimInVecParam} \times \mathbb{R}^{\DimSOSBasis\times \DimSOSBasis} \times  \mathbb{R}^{\DimASparam\times \NumAffine}$,
	there exists a unique solution ${\SOSdotVMat{\IDAffine}}(\InVecParam, \LyapMat,\ASparam)$ of ${\SOSdotVMat{\IDAffine}}$ that satisfies \eqref{eq:def_dotLyapF_new}--\eqref{eq:def_symmetricSOSdotVMat}.
Furthermore, this solution is bilinear in $(\InVecParam,\LyapMat)$ and linear in $\ASparam$.

\end{lemma}
\begin{\MyProof}
	The proof is described in Appendix \ref{pf:design_penaF}.
\end{\MyProof}
\begin{remark}[{\MyHighlight{Contributions of  Lemma \ref{thm:design_penaF}}}]
Owing to the uniqueness of ${\SOSdotVMat{\IDAffine}}(\InVecParam, \LyapMat,\ASparam)$ from Lemma \ref{thm:design_penaF}, 
the penalty function $\PenaltyFuncRegu(\InVecParam,\LyapMat,\ASparam)$ in \eqref{eq:def_PenaltyFuncRegu} is well-defined.	
It is straightforward to determine the constants $(\AScoefmat,\AScoefscalar)$ because  \eqref{eq:def_dotLyapF_new}--\eqref{eq:def_symmetricSOSdotVMat}  reduce to simple linear equations of ${\SOSdotVMat{\IDAffine}}$ independent of $\sysx$, as demonstrated in Section \ref{sec_plant}.
\end{remark}

Using the well-defined penalty function $\PenaltyFuncRegu(\InVecParam,\LyapMat,\ASparam)$, we derive the following results: 
The gradient method guarantees the stability.

\begin{theorem}[{\MyHighlight{Proper gradient method}}]\label{thm:stabilizing_optim}
	Suppose that Assumptions \ref{ass:polypoly_sys} and \ref{ass:SOS} hold.
	Given 
	initial decision variables $({\Ite{\InVecParam}{0}},{\Ite{\LyapMat}{0}},{\Ite{\ASparam}{0}})$
	and
	constants $\iniPenaltyFuncRegu \in \mathbb{R}$ in \eqref{eq:def_PenaltyFuncRegu} 
	and
	$\appcostLB \in \mathbb{R}$,
	suppose that the following conditions hold:
	\begin{align}
	{\Ite{\LyapMat}{0}} + {\Ite{\LyapMat}{0}}^{\Mytop} &\succ 0
	, \label{eq:cond_robust_stab_P}
	\\
	{\SOSdotVMat{\IDAffine}}({\Ite{\InVecParam}{0}},{\Ite{\LyapMat}{0}},{\Ite{\ASparam}{0}}) &\succ 0
	, \label{eq:cond_robust_stab_T}
	\\
	\ObjF(  {\Ite{\InVecParam}{0}}   ,{\Ite{\LyapMat}{0}} ,{\Ite{\ASparam}{0}}) 
	&<
	\appcostLB
	+ \iniPenaltyFuncRegu
	.  \label{eq:cond_iniPenaltyFuncRegu}
	\end{align}
	For any $\IDIteGrad\in\{0,1,2,\dots\} $, if the gradient method in \eqref{eq:gradient_updatate_InVecParam}--\eqref{eq:gradient_updatate_ASparam} provides $(  {\Ite{\InVecParam}{\IDIteGrad}}   ,{\Ite{\LyapMat}{\IDIteGrad}} ,{\Ite{\ASparam}{\IDIteGrad}})$ satisfying
	\begin{align}
	\ObjF(  {\Ite{\InVecParam}{\IDIteGrad}}   ,{\Ite{\LyapMat}{\IDIteGrad}} ,{\Ite{\ASparam}{\IDIteGrad}})
	&\leq
	\ObjF(  {\Ite{\InVecParam}{0}}   ,{\Ite{\LyapMat}{0}} ,{\Ite{\ASparam}{0}})
	,\label{pf:ObjF_mono_dec}
	\\
	\Expect{\StoParam}[\Expect{\IniState}[    
	\EachValue(\IniState,\VParam(  {\Ite{\InVecParam}{\IDIteGrad}}    ,\StoParam))
	] ] & \geq \appcostLB
	,  \label{eq:cost_morethan_LB}
	\end{align}
	then the robust global asymptotic stability of the feedback system \eqref{eq:def_sys} with $\Input(\sysx)=\InpOpt(\sysx, {\Ite{\InVecParam}{\IDIteGrad}})$ is satisfied:
	\begin{align}
	\forall 
	\StoParam  \in \SetSto
	,\;
	\forall 
	\IniState  \in \mathbb{R}^{\DimX}
	,\;
	\lim_{\MyT \rightarrow \infty } \sysx(\MyT,\StoParam) = 0 
	. \label{eq:_def_robust_global_asymptotic_stable}
	\end{align}

\end{theorem}
\begin{\MyProof}
	The proof is described in Appendix \ref{pf:stabilizing_optim}.
\end{\MyProof}
\begin{remark}[{\MyHighlight{Contributions of  Theorem \ref{thm:stabilizing_optim}}}]
	Theorem \ref{thm:stabilizing_optim} implies that the designed controller satisfies the robust global asymptotic stability \eqref{eq:_def_robust_global_asymptotic_stable}, while the gradient method optimizes the controller parameter $\InVecParam$.
	The gradient method can satisfy the condition \eqref{pf:ObjF_mono_dec} naturally.
	The conditions \eqref{eq:cond_iniPenaltyFuncRegu} and \eqref{eq:cost_morethan_LB} are easy to satisfy by setting $\appcostLB$ and $\iniPenaltyFuncRegu$ as a sufficiently small value and large value, respectively. 
	Because Theorem \ref{thm:stabilizing_optim} requires initial decision variables $({\Ite{\InVecParam}{0}},{\Ite{\LyapMat}{0}},{\Ite{\ASparam}{0}})$ to satisfy \eqref{eq:cond_robust_stab_P} and \eqref{eq:cond_robust_stab_T},
	they are designed in the following.	
\end{remark}

We design the initial decision variables as solutions to the following two optimization problems \eqref{eq:SDP_to_LyapInv_tmpInpParam} and \eqref{eq:SDP_to_ASparam}:
\begin{align}
	&
	\max_{   \LyapInv, \tmpInpParam, {\SOSRtdLyapMat{1}}, \dots, {\SOSRtdLyapMat{\NumAffine}} , \SOSAepsilon }
	\SOSAepsilon
	\nonumber\\&
	\mathrm{s.t.}
	\begin{cases}
		\forall \IDAffine,\;\;
		\forall \sysx,\;\;
		{\tdLyapMat{\IDAffine}}(\sysx,\LyapInv, \tmpInpParam)
		+{\tdLyapMat{\IDAffine}}(\sysx,\LyapInv, \tmpInpParam)^{\Mytop}
		\\
		\qquad\qquad\;\;
		=
		( \SOSmatBasis(\sysx)  \otimes  \Identity{\DimSOSBasis} )^{\Mytop}
		{\SOSRtdLyapMat{\IDAffine}}
		( \SOSmatBasis(\sysx)  \otimes  \Identity{\DimSOSBasis} )
		,\\
		\forall \IDAffine,\;\;
		{\SOSRtdLyapMat{\IDAffine}} = {\SOSRtdLyapMat{\IDAffine}^{\Mytop}} 
		\succeq \SOSAepsilon {\Identity{\DimSOSmatBasis\DimSOSBasis}}
		,\\		
		\LyapInv   = \LyapInv^{\Mytop}
		\succeq  \SOSAepsilon {\Identity{\DimSOSBasis}}
		,\\
		\SOSAepsilon > 0
		.
	\end{cases}
	\label{eq:SDP_to_LyapInv_tmpInpParam}
\end{align}
where
\begin{align}
{\tdLyapMat{\IDAffine}}(\sysx,\LyapInv, \tmpInpParam)
&:=
-
\frac{\partial \SOSBasis(\sysx)}{\partial \sysx^{\Mytop} }
\Big(
{\PolyDriftMat{\IDAffine}}(\sysx) 
\LyapInv
+
{\PolyInMat{\IDAffine}}(\sysx) 
\InBasis(\sysx)
\tmpInpParam	
\Big)
.\label{eq:def_tdLyapMat}
\end{align} 
If a feasible solution  $(\solLyapInv, \soltmpInpParam )$
to the problem \eqref{eq:SDP_to_LyapInv_tmpInpParam} is obtained, 
the initial decision variables are defined as follows: 
\begin{align}
{\Ite{\InVecParam}{0}} &:={\VEC{  \soltmpInpParam \solLyapInv^{-1}  }}
, \label{eq:def_ini_InVecParam}
\\
{\Ite{\LyapMat}{0}} &:=\solLyapInv^{-1}
, \label{eq:def_ini_LyapMat}
\\
{\Ite{\ASparam}{0}} &:=\solASparam
, \label{eq:def_ini_ASparam}
\end{align} 
where $\solASparam$ is a feasible solution to the following problem:
\begin{align}
\max_{\ASparam,{\dotVepsilon}}
{\dotVepsilon}
\;\;
\mathrm{s.t.}
\begin{cases}
\forall \IDAffine,\;\;
{\SOSdotVMat{\IDAffine}}({\Ite{\InVecParam}{0}}, {\Ite{\LyapMat}{0}},\ASparam)
\succeq {\dotVepsilon} {\Identity{\DimSOSdotVBasis}}
,\\
{\dotVepsilon} > 0
,
\end{cases}
\label{eq:SDP_to_ASparam}
\end{align}
Note that ${\SOSdotVMat{\IDAffine}}({\Ite{\InVecParam}{0}}, {\Ite{\LyapMat}{0}},\ASparam)$ is the unique solution satisfying \eqref{eq:def_dotLyapF_new}--\eqref{eq:def_symmetricSOSdotVMat} for 
$(\InVecParam, \LyapMat)=({\Ite{\InVecParam}{0}}, {\Ite{\LyapMat}{0}})$.
We derive the following theorem to guarantee the feasibility of these optimization problems and the reasonableness of the initial decision variables.

\begin{theorem}[{\MyHighlight{Initial controller design}}]\label{thm:stabilizing_initial}
Supposing that Assumptions \ref{ass:polypoly_sys} and \ref{ass:SOS} hold, the following properties are satisfied:
\begin{enumerate}
	\item 
	[{\MyLabelFeasibleWP}]
	There exists a feasible solution $(\solLyapInv, \soltmpInpParam, {\solSOSRtdLyapMat{1}}, \dots, {\solSOSRtdLyapMat{\NumAffine}} ,\solSOSAepsilon)$
	to the problem \eqref{eq:SDP_to_LyapInv_tmpInpParam}.
	
	\item
	[{\MyLabelFeasibleAddP}]
	Given a feasible solution to the problem \eqref{eq:SDP_to_LyapInv_tmpInpParam},
	there exists a feasible solution  $\solASparam$ to the problem \eqref{eq:SDP_to_ASparam}
	for some $\AScoefmat \in \mathbb{R}^{\DimASparam \times \DimSOSdotVBasis(\DimSOSdotVBasis+1)/2 } $ and $\AScoefscalar\in \{0,1\}$ in \eqref{eq:def_ASconstraint}.
	
	\item
	[{\MyLabelInitialVal}]
	The initial decision variables $({\Ite{\InVecParam}{0}}, {\Ite{\LyapMat}{0}}, {\Ite{\ASparam}{0}})$ in
	\eqref{eq:def_ini_InVecParam}--\eqref{eq:def_ini_ASparam}
	satisfy the conditions \eqref{eq:cond_robust_stab_P} and \eqref{eq:cond_robust_stab_T}.

\end{enumerate}

\end{theorem}
\begin{\MyProof}
	The proof is described in Appendix \ref{pf:stabilizing_initial}.		
\end{\MyProof}
\begin{remark}[{\MyHighlight{Contributions of  Theorem \ref{thm:stabilizing_initial}}}]
Theorem \ref{thm:stabilizing_initial} ensures that the initial decision variables defined in \eqref{eq:def_ini_InVecParam}--\eqref{eq:def_ini_ASparam} are always obtained and satisfy the conditions \eqref{eq:cond_robust_stab_P} and \eqref{eq:cond_robust_stab_T} required by Theorem \ref{thm:stabilizing_optim}.
Consequently, the proposed gradient method yields stabilizing suboptimal controllers. 
Moreover, the optimization problems \eqref{eq:SDP_to_LyapInv_tmpInpParam} and \eqref{eq:SDP_to_ASparam} are easy to solve because they belong to a class of convex optimization problems.
Specifically, the SOS constraints can be expressed as linear equations of the decision variables, reducing \eqref{eq:SDP_to_LyapInv_tmpInpParam} and \eqref{eq:SDP_to_ASparam} to convex semidefinite programming (SDP).
\end{remark}

\begin{remark}[{\MyHighlight{Comparison with related work}}]\label{rem:related_work}
	While our strategy employing SOS in \eqref{eq:def_Lyap} and \eqref{eq:def_dLyap_dt} builds on groundbreaking results \citep{Xu09TAC,Prajna04}, it offers the novelties and advantages as presented above:
	The proposed method guarantees the robust stability of the system regarding the unknown stochastic parameter $\StoParam$ whereas 
	\citep[Theorem 1]{Xu09TAC} and \citep[Theorem 6]{Prajna04} do not treat $\StoParam$.
	In addition, the proposed penalty function $\PenaltyFuncRegu(\InVecParam,\LyapMat,\ASparam)$ enables the integration of the stability guarantee with the minimization of the expected cost function via the gradient method. 
	In contrast, \citep{Xu09TAC,Prajna04} focus solely on stability without addressing cost minimization.
\end{remark}

\subsection{Implementation} \label{sec_implementation}

\begin{algorithm}[!t]  
	\renewcommand{\algorithmicrequire}{\textbf{Input:}}
	\renewcommand{\algorithmicensure}{\textbf{Output:}}                       
	\caption{Design of the stochastic suboptimal feedback controller}  
	\label{alg1}
	\begin{algorithmic}[1]  
			
		\Require $\xCost(\sysx)$, $\uCost(\sysx)$, $\PDF{\StoParam}$, $\PDF{\IniState}$, $\JBasis(\sysx)$,
		$\MarginalWeight(\sysx)$, $\CoefReguVParam$, $\penaW$, $\iniPenaltyFuncRegu$, and $\NumIteGrad$ 
		
		\Ensure Suboptimal controller $\InpOpt(\sysx, \Ite{\InVecParam}{\NumIteGrad})$

		
		\noindent
		\textbf{--Design of the initial decision variables--}	
		
		\State
		Define 
		$\SOSBasis(\sysx)$,
		$\SOSmatBasis(\sysx)$,
		$\SOSdotVBasis(\sysx)$, and 
		$\InBasis(\sysx)$
		in Remark \ref{rem:monomial_bases} according to Assumption \ref{ass:SOS}.
		
		\State
		Obtain ${\Ite{\InVecParam}{0}}$ and ${\Ite{\LyapMat}{0}}$ from \eqref{eq:def_ini_InVecParam} and \eqref{eq:def_ini_LyapMat} by solving the SDP \eqref{eq:SDP_to_LyapInv_tmpInpParam}
		
		\State
		Obtain ${\Ite{\ASparam}{0}}$ from \eqref{eq:def_ini_ASparam} by solving the SDP \eqref{eq:SDP_to_ASparam}
		
		\noindent
		\textbf{--Optimization of the controller with guaranteeing the stability--}	

		\State
		Calculate $\Expect{\IniState}[ \JBasis(\IniState) ]$

		\State
		Calculate the constant matrices 
		$\MarginalMap[]{ 
			{\xAterm{\IDAffine}}
			{\xAterm{\IDbAffine}}^{\Mytop}
		}$
		and
		$\MarginalMap[]{ 
			{\xAterm{\IDAffine}}
			{\xBterm}^{\Mytop}
		}$
		according to Theorem \ref{thm:explicit_cost}

		\For{$\IDIteGrad=0$ to $\NumIteGrad-1$}

		\State \label{algLINE:grad}
		Calculate the gradient of $\ObjF(  {\Ite{\InVecParam}{\IDIteGrad}}   ,{\Ite{\LyapMat}{\IDIteGrad}} ,{\Ite{\ASparam}{\IDIteGrad}})$ in \eqref{eq:def_ObjF} by  
		Theorem \ref{thm:explicit_grad} and Remarks \ref{eq:approx_expectation} and \ref{eq:penaF_derivative}

		\State \label{algLINE:StepSize}
		Obtain the step size $ {\CoefGrad{\IDIteGrad}}$ according to Remark \ref{rem:stepsize}
		
		\State
		Update the parameters $(  {\Ite{\InVecParam}{\IDIteGrad+1}}   ,{\Ite{\LyapMat}{\IDIteGrad+1}} ,{\Ite{\ASparam}{\IDIteGrad+1}})$ by \eqref{eq:gradient_updatate_InVecParam}--\eqref{eq:gradient_updatate_ASparam}
		
		\EndFor

	\end{algorithmic}
\end{algorithm}

Lastly, the proposed method is summarized in Algorithm \ref{alg1}, which consists of two main parts: the design of the initial decision variables and the optimization of the controller.
The following remarks elaborate on several technical aspects of the algorithm.

\begin{remark}[{\MyHighlight{Approximation of the expectations}}]\label{eq:approx_expectation}
	We aim to approximate the expectation $\Expect{\StoParam}[\dots]$ used in the gradient method in \eqref{eq:gradient_updatate_InVecParam}
	because the expectations for many types of PDFs $\PDF{\StoParam}$ are infeasible to calculate in an exact sense.
	The expectations can be approximated using several methods such as a numerical integration with the trapezoidal rule and the Monte Carlo approximation. 
	Stochastic gradient descent methods \citep[Sec. III]{Sun20Cyb},\citep{Adam-Kingma15,JMLR:v18:17-049} are efficient for approximating the expected gradient.
	Using these techniques approximates the gradient in \eqref{eq:gradient_updatate_InVecParam} by
	${\partial }
	\Expect{\StoParam}[\Expect{\IniState}[    
	\EachValue(\IniState,\VParam(    {\Ite{\InVecParam}{\IDIteGrad}}   ,\StoParam))
	]  ] 
	/{\partial \InVecParam } 
	\approx
	({1}/{\NumMC})
	\sum_{\IDMC=1}^{\NumMC}
	{\partial }
	\Expect{\IniState}[    
	\EachValue(\IniState,\VParam(    {\Ite{\InVecParam}{\IDIteGrad}}   ,{\MC{\StoParam}{\IDMC,\IDIteGrad}}))
	]
	/{\partial \InVecParam } 
	$,
	where ${\MC{\StoParam}{\IDMC,\IDIteGrad}}$ are parameters sampled according to the approximation techniques.
\end{remark}

\begin{remark}[{\MyHighlight{Gradient of the penalty function}}]\label{eq:penaF_derivative}
	In Line \ref{algLINE:grad}, the gradient of the penalty function is obtained using the following relation \citep[Eq. (C.22)]{Bishop06}: 
	\begin{align}
	\frac{\partial}{\partial {\El{\LyapMat}{\IDEl,\IDbEl}}}
	\Big(
	-	\ln \det(\LyapMat)	
	\Big)	
	=
	- {\TRACE{
			\LyapMat^{-1}  \frac{\partial \LyapMat}{\partial {\El{\LyapMat}{\IDEl,\IDbEl}}}
	}}
	.		
	\end{align}
	The gradient of $\det({\SOSdotVMat{\IDAffine}}(\InVecParam, \LyapMat,\ASparam) )	$ in the penalty function is obtained in the same manner.
\end{remark}

\begin{remark}[{\MyHighlight{Step size}}]\label{rem:stepsize}
In Line \ref{algLINE:StepSize}, it is advisable to determine the step size ${\CoefGrad{\IDIteGrad}}$ used in \eqref{eq:gradient_updatate_InVecParam}--\eqref{eq:gradient_updatate_ASparam} such that the condition \eqref{pf:ObjF_mono_dec} is satisfied.
A suitable approach is to employ a step size ${\CoefGrad{\IDIteGrad}}$ that obeys the Wolfe condition \citep[Section 3.1]{Nocedal06}:
\begin{align}
\WolfeObjF({\Ite{\WolfeVar}{\IDIteGrad+1}})
&\leq
\WolfeObjF({\Ite{\WolfeVar}{\IDIteGrad}})
-
\WolfeAcoef
{\CoefGrad{\IDIteGrad}}
\Big\|
\frac{\partial \WolfeObjF}{\partial \WolfeVar}
({\Ite{\WolfeVar}{\IDIteGrad}})
\Big\|^{2}		
,\label{eq:Wolfe_A}
\end{align}	
where $\WolfeVar$ and $\WolfeObjF({\Ite{\WolfeVar}{\IDIteGrad}})$ denote the vector-valued collection of $(\InVecParam, \LyapMat,\ASparam)$ and  the objective function $\ObjF(  {\Ite{\InVecParam}{\IDIteGrad}}   ,{\Ite{\LyapMat}{\IDIteGrad}} ,{\Ite{\ASparam}{\IDIteGrad}})$, respectively, for the brevity of the notation. 
The symbol $\WolfeAcoef \in (0,1)$ is a free parameter.
A backtracking approach \citep[Algorithm 3.1]{Nocedal06} can be employed to find such a step size ${\CoefGrad{\IDIteGrad}}$.
At each iteration $\IDIteGrad$, we firstly set ${\CoefGrad{\IDIteGrad}}=\WolfeIniStepSize$, where $\WolfeIniStepSize$ is a large positive value. We next repeat ${\CoefGrad{\IDIteGrad}} \leftarrow \WolfeDiscount {\CoefGrad{\IDIteGrad}}$ until \eqref{eq:Wolfe_A} holds, where $\WolfeDiscount \in (0,1)$ is a free parameter.
\end{remark}

\section{Numerical example} \label{sec_numerical_example}

This section evaluates the effectiveness of the proposed method.
Section \ref{sec_plant} introduces the target system and outlines the simulation settings. 
Subsequently, Section \ref{sec_results} presents and analyses the simulation results.

\subsection{Plant system and setting}\label{sec_plant}

Let us consider the following example of the system \eqref{eq:def_sys} with $\DimX=2$:
\begin{align}	
\Drift(\sysx,\StoParam)
&:=
\begin{bmatrix}
{\El{\sysx}{1}} + {\El{\sysx}{1}^{2}} + ({\El{\StoParam}{1}}-2) {\El{\sysx}{1}^{3}} - {\El{\sysx}{1}}{\El{\sysx}{2}^{2}}/2
+
{\El{\sysx}{2}}
\\	
({\El{\StoParam}{2}}+1) {\El{\sysx}{1}}+ {\El{\sysx}{2}^{2}}
\end{bmatrix}
,\label{eq:sim_Drift}
\\
\InMat(\StoParam)
&:=
\begin{bmatrix}
{\El{\StoParam}{1}}\\ 
{\El{\StoParam}{2}}+1 
\end{bmatrix}
,
\end{align}
where the uncertain parameter $\StoParam$ obeys the uniform distribution on the two dimensional finite set 
$\SetSto:=\{ 0,0.1,0.9,1.0  \}^{2} $.
This example satisfies Assumption \ref{ass:polypoly_sys} with $\NumAffine=4$ by the following setting:
\begin{align*}
\forall \IDAffine \in\{1,2,3,4\},\quad
\xDrift{\IDAffine}(\sysx)&=\Drift(\sysx,{\vertexStoParam{\IDAffine}})
,\\
\forall \IDAffine \in\{1,2,3,4\},\quad
\xInMat{\IDAffine} &=\InMat({\vertexStoParam{\IDAffine}})
,\\
{\sDrift{1}}(\StoParam)&=(1-{\El{\StoParam}{1}})(1-{\El{\StoParam}{2}})
,\\
{\sDrift{2}}(\StoParam)&=(1-{\El{\StoParam}{1}}){\El{\StoParam}{2}}
,\\
{\sDrift{3}}(\StoParam)&={\El{\StoParam}{1}}(1-{\El{\StoParam}{2}})
,\\
{\sDrift{4}}(\StoParam)&={\El{\StoParam}{1}}{\El{\StoParam}{2}}
,
\end{align*}
where 
${\vertexStoParam{1}}:=[0,0]^{\Mytop}$,
${\vertexStoParam{2}}:=[0,1]^{\Mytop}$,
${\vertexStoParam{3}}:=[1,0]^{\Mytop}$, and
${\vertexStoParam{4}}:=[1,1]^{\Mytop}$.

We define the monomial bases
$\SOSBasis(\sysx)$,
$\SOSmatBasis(\sysx)$,
$\InBasis(\sysx)$, and 
$\SOSdotVBasis(\sysx)$
in Remark \ref{rem:monomial_bases} according to Assumption \ref{ass:SOS} as follows:
\begin{align}
\SOSBasis(\sysx)
&:=[
{\El{\sysx}{1}},
{\El{\sysx}{2}}
]^{\Mytop}
,\\
\SOSmatBasis(\sysx)
&:=
[
1,
{\El{\sysx}{1}},
{\El{\sysx}{2}}
]^{\Mytop}
,\\
\InBasis(\sysx) 
&:= 
[
1,
{\El{\sysx}{1}},
{\El{\sysx}{2}},
{\El{\sysx}{1}^{2}},
{\El{\sysx}{1}}{\El{\sysx}{2}},
{\El{\sysx}{2}^{2}}
]
,\\
\SOSdotVBasis(\sysx) 
&:= 
[
{\El{\sysx}{1}},
{\El{\sysx}{2}},
{\El{\sysx}{1}^{2}},
{\El{\sysx}{1}}{\El{\sysx}{2}},
{\El{\sysx}{2}^{2}}
]^{\Mytop}
.
\end{align}
Because $\Drift(\sysx,\StoParam)$ in \eqref{eq:sim_Drift} is decomposed into $\Drift(\sysx,\StoParam) = \DriftMat(\sysx,\StoParam) \SOSBasis(\sysx)$ with
\begin{align}
\DriftMat(\sysx,\StoParam)
&:=
\begin{bmatrix}
1+ {\El{\sysx}{1}} + ({\El{\StoParam}{1}}-2) {\El{\sysx}{1}^{2}}  - {\El{\sysx}{2}^{2}}/2  \;&\; 1 \\
{\El{\StoParam}{2}}+1 \;&\; {\El{\sysx}{2}}
\end{bmatrix}
,
\end{align}	
the matrix ${\PolyDriftMat{\IDAffine}}(\sysx) $ satisfying $ \xDrift{\IDAffine}(\sysx) = {\PolyDriftMat{\IDAffine}}(\sysx) \SOSBasis(\sysx)$ in \eqref{eq:def_PolyDriftMat} is given by ${\PolyDriftMat{\IDAffine}}(\sysx)=\DriftMat(\sysx,{\vertexStoParam{\IDAffine}})$.
Under these settings, the constraint \eqref{eq:def_ASconstraint} reduces to the three equations:
${\El{\SOSdotVMat{\IDAffine}}{5,3}}
= {\El{\ASeach{\IDAffine}}{1}}$,
${\El{\SOSdotVMat{\IDAffine}}{3,2}}
= {\El{\ASeach{\IDAffine}}{2}}$, and
${\El{\SOSdotVMat{\IDAffine}}{4,2}}
= {\El{\ASeach{\IDAffine}}{3}}$.

In the design of the stochastic optimal controller,
the cost function in \eqref{eq:def_cost} is defined with $\xCost( \sysx )=\sysx^{\Mytop}\sysx$ and $\uCost=10$.
The PDF $\PDF{\IniState}$ is set as the uniform distribution on $\{ -3,0,3 \}^{2} \setminus \{[0,0]^{\Mytop}\}$.
The basis $\JBasis(\sysx)$ consists of all the monomials whose degrees are greater than $0$ and less than or equal to ${6}$, that is, ${\El{\sysx}{1}^{\simBasisDegA}}{\El{\sysx}{2}^{\simBasisDegB}}$ for all $(\simBasisDegA,\simBasisDegB) \in \{0,1,2,\dots,6\}^{2}$ satisfying $0 < \simBasisDegA + \simBasisDegB \leq {6}$. 
The coefficient in \eqref{eq:def_OptVParam} is set to $\CoefReguVParam=0$. 
The weight function in \eqref{eq:def_MarginalMap} is defined by $\MarginalWeight(\sysx):= \sum_{\IDProj=1}^{\NumProj} \MyDeltaFunc(\sysx-\DeltaSampleState{\IDProj})$, where $\MyDeltaFunc(\NonArg)$ is the Dirac delta function and
$\DeltaSampleState{\IDProj}$ are set to all the members of the set $\{-3,-2.9,-2.8,\dots,3\}^{2}$.
In the SDP \eqref{eq:SDP_to_LyapInv_tmpInpParam}, the SOS constraints are replaced with linear equations ${\SOSalgebraicEq{\IDAffine}}(\LyapInv, \tmpInpParam, {\SOSRtdLyapMat{\IDAffine}})=0$, where they are numerically implemented as $ |{\El{\SOSalgebraicEq{\IDAffine}(\LyapInv, \tmpInpParam, {\SOSRtdLyapMat{\IDAffine}})}{\IDEl}} |\leq 1.0 \times 10^{-15}$.

In the gradient method,
the number of the iterations is set to $\NumIteGrad={2000}$.
The constants used in the penalty function in \eqref{eq:def_PenaltyFuncRegu} are set to $\penaW=0.1$ and $\iniPenaltyFuncRegu=1 \times 10^{20}$. 
To determine the step size ${\CoefGrad{\IDIteGrad}}$ according to the condition \eqref{eq:Wolfe_A}, we set the parameters to $\WolfeAcoef=1.0\times 10^{-4}$, $\WolfeIniStepSize=0.01$, and $\WolfeDiscount=0.5$.

\subsection{Simulation results} \label{sec_results}

\newcommand{\simFolder}{results_Art3_20241218}

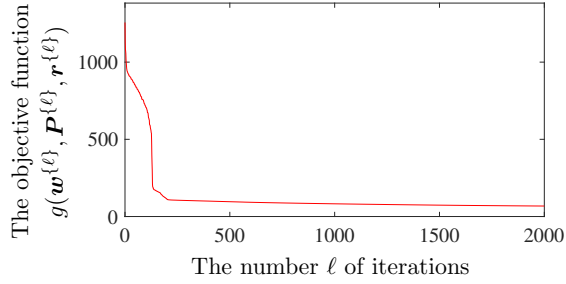
\begin{figure}[!t]	
	\vspace*{+0.07in}
	\centering
	\begin{tikzpicture}[scale=0.9]%
		\begin{axis}[axis y line=none, axis x line=none
			,xtick=\empty, ytick=\empty 
			,xmin=-10,xmax=10,ymin=-10,ymax=10
			,width=4.0 in,height=2.3 in
			]	
			
			
			\node at (0.5,0) { 
				\includegraphics[width=0.49\linewidth]{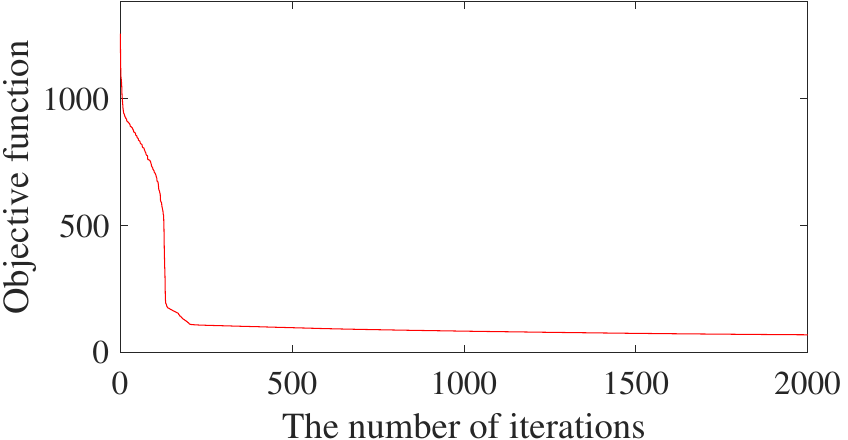} 
			};
			\fill[fill=white] (-11,-10) rectangle (-7.5,10);
			\node at (-8.8,1.0) {\rotatebox{90}{
				\shortstack{ The objective function \\ 
					$\ObjF(  {\Ite{\InVecParam}{\IDIteGrad}}   ,{\Ite{\LyapMat}{\IDIteGrad}} ,{\Ite{\ASparam}{\IDIteGrad}})$ }
			}};
			\fill[fill=white] (-6,-7.9) rectangle (6,-11);
			\node at (1.2,-9) {The number $\IDIteGrad$ of iterations};		
		\end{axis}			
	\end{tikzpicture}	
	
	\caption{Results of the optimization \eqref{eq:optim_param} by using the proposed gradient method}
	\label{fig:optim_transition}
\end{figure}

\newcommand{\FigTemplate}[1]{%
	\begin{tikzpicture}[scale=0.9]%
	\begin{axis}[axis y line=none, axis x line=none
	,xtick=\empty, ytick=\empty 
	,xmin=-10,xmax=10,ymin=-10,ymax=10
	,width=2.4 in,height=2.2 in
	]	
	
	
	\node at (-0.0,-0.2) { 
		\includegraphics[width=1.05\linewidth]{#1} 
	};
	\fill[fill=white] (-11,-6) rectangle (-8.0,6);
	\node at (-8.8,1.0) {\rotatebox{90}{
			Second state $\El{\sysx(\MyT,\StoParam)}{2}$
	}};
	\fill[fill=white] (-6,-7.9) rectangle (6,-11);
	\node at (1.2,-9) {First state $\El{\sysx(\MyT,\StoParam)}{1}$};		
	\end{axis}			
	\end{tikzpicture}
	\vspace*{-0.3\baselineskip}	
}

\newcommand{\MyBBsize}{0.27\linewidth}
\begin{figure}[!t]	
	\vspace*{+0.07in}
	\centering
	\begin{minipage}[b]{\MyBBsize}
		\centering	
		{\FigTemplate{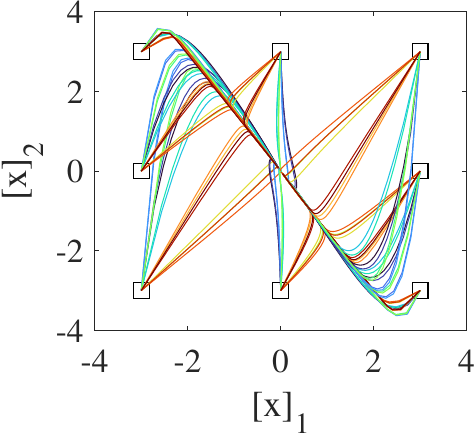}}		
		\subcaption{Without optimization}
		\label{fig:plot_wo_optimization}
	\end{minipage}%
	\begin{minipage}[b]{\MyBBsize}
		\centering	
		{\FigTemplate{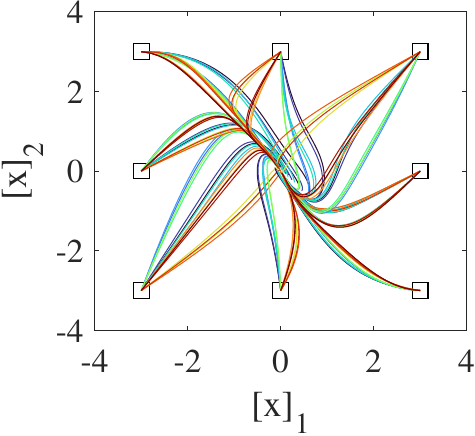}}
		\subcaption{Without optimality}
		\label{fig:plot_wo_optimality}	
	\end{minipage}%
	
	\begin{minipage}[b]{\MyBBsize}
		\centering	
		{\FigTemplate{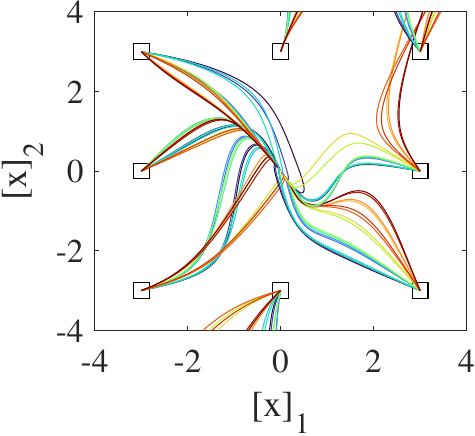}}
		\subcaption{Without stability}
		\label{fig:plot_wo_stability}	
	\end{minipage}%
	\begin{minipage}[b]{\MyBBsize}
		\centering	
		{\FigTemplate{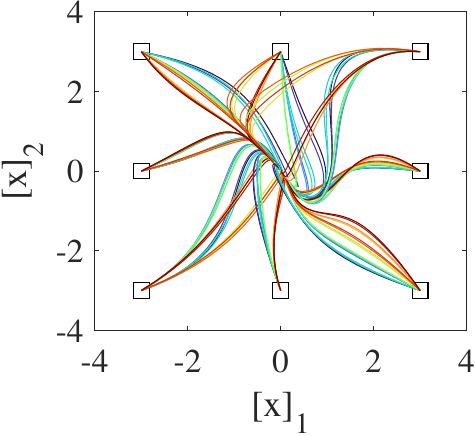}}
		\subcaption{Proposed controller}
		\label{fig:plot_proposed}	
	\end{minipage}%
	
	\caption{Control results of the proposed controller compared with the three controllers.
	The square markers denote the initial states $\sysx(0,\StoParam)$.
	The different colored trajectories express various values of the uncertain parameter $\StoParam$.}
	\label{fig:state_transition}
\end{figure}

Figure \ref{fig:optim_transition} illustrates the results of the optimization process \eqref{eq:optim_param} for designing the proposed controller.
The objective function $\ObjF(\InVecParam,\LyapMat,\ASparam)$ was successfully decreased using the gradient method, resulting in the stabilizing suboptimal controller  $\InpOpt(\sysx, \Ite{\InVecParam}{\NumIteGrad})$.

The performance of the proposed stabilizing suboptimal controller $\InpOpt(\sysx, \Ite{\InVecParam}{\NumIteGrad})$ is compared with three baseline controllers.
Firstly, a controller \textit{without optimization} corresponds to the stabilizing controller  $\InpOpt(\sysx, {\Ite{\InVecParam}{0}})$ where the gradient-based optimization in \eqref{eq:gradient_updatate_InVecParam}--\eqref{eq:gradient_updatate_ASparam} is not performed. 
Secondly, a controller \textit{without optimality} is designed by replacing the objective function in \eqref{eq:optim_param} with $\|\InVecParam\|^{2} + \PenaltyFuncRegu(\InVecParam, \LyapMat,\ASparam)$, relying on  the quadratic controller gain $\|\InVecParam\|^{2}$ instead of the approximate cost function $\EachValue(\IniState,\VParam(\InVecParam,\StoParam))$. 
Thirdly, a controller \textit{without stability} is obtained by removing the penalty function $\PenaltyFuncRegu(\InVecParam, \LyapMat,\ASparam)$ 
from the objective function, where the gradient method is terminated when the approximate cost becomes negative.

In Fig. \ref{fig:state_transition}, 
\subref{fig:plot_wo_optimization}, \subref{fig:plot_wo_optimality}, \subref{fig:plot_wo_stability}, and \subref{fig:plot_proposed} correspond to the controllers without optimization, without optimality, without stability, and the proposed stabilizing suboptimal controller, respectively. 
The controllers in \subref{fig:plot_wo_optimization}, \subref{fig:plot_wo_optimality}, and \subref{fig:plot_proposed} successfully stabilized the feedback system for all uncertain parameters $\StoParam$ and initial states $\IniState$.
These results confirm the effectiveness of the proposed penalty function $\PenaltyFuncRegu(\InVecParam, \LyapMat,\ASparam)$.
Table \ref{tab:J_results} summarizes the control performance in terms of the expected cost function $\Expect{\StoParam}[
\Expect{\IniState}[    
\costJ(\Input,\IniState,\StoParam) 
] 
]$, where the infinite time horizon is approximated as the finite time ${\terminalT=30}$ in the simulation. 
The proposed controller achieves a lower expected cost compared to the baseline controllers, confirming that the proposed gradient method with parametric approximation improves control performance.

\begin{table}
	\centering
	\caption{Expected cost function $\Expect{\StoParam}[\Expect{\IniState}[ \costJ(\Input, \sysx,\StoParam) ] ]$ indicating the average control performance regarding $\StoParam$ and $\IniState$.}
	\label{tab:J_results}
	{\scriptsize
	\begin{tabular}{| c|c |}
		\hline
		Controller &  Expected cost   \\		\hline\hline
		Without optimization	 &  1262 \\ 	\hline	 
		Without optimality	 &  157 \\ 	\hline	 
		Without stability	 &  diverged \\ 	\hline	 
		\textbf{Stabilizing suboptimal (Proposed)}	 & \textbf{74} \\ 	\hline		
	\end{tabular}
	}
\end{table}

\section{Conclusion} \label{sec_conclusion}

This study presented a method for designing suboptimal controllers for time-invariant nonlinear stochastic systems, satisfying three essential requirements: feedback control, offline design, and robust stabilization.
We overcame the challenges associated with nonlinear optimal control, time-invariant stochastic characteristics, and robust stability by proposing a parametric approximation of cost functions, gradient-based controller optimization, and an SOS-based penalty function, respectively.
The proposed design method is supported by key Theorems \ref{thm:BellmanResidual}--\ref{thm:stabilizing_initial}.
A residual in approximating a cost function is analyzed in Theorem \ref{thm:BellmanResidual}.
An optimal approximate cost function and its gradient are explicitly derived by Theorem \ref{thm:explicit_cost} and \ref{thm:explicit_grad}, which enables us to employ the gradient method.
Theorem \ref{thm:stabilizing_optim} shows that the robust stability is guaranteed while the gradient method is employed.
Appropriate initial decision variables used in the gradient method are obtained using Theorem \ref{thm:stabilizing_initial}.

Future research will focus on extending the proposed method to other classes of nonlinear stochastic systems, such as
distributionally uncertain nonlinear stochastic systems, a linear case of which has been analyzed in \citep{Ito24TAC}.
Other promising directions include extensions for alternative control problems, such as risk-sensitive control \citep{ItoTAC19}, and applications to practical systems, including nanoscale devices with manufacturing variations \citep{ItoACCESS19}.



\appendix
\section*{Appendix}
\addcontentsline{toc}{section}{Appendix}

\section{Proof of Theorem \ref{thm:BellmanResidual}}\label{pf:BellmanResidual}
For the brevity of the notation, let us define
$
\BellA(\sysx):=
\xCost(\sysx)
+
\InpOpt(\sysx, \InVecParam)^{\Mytop}
\uCost(\sysx)
\InpOpt(\sysx, \InVecParam)/2
\geq 0$
and
$
\BellB(\sysx):=
(
{\partial \EachValue(\sysx,\preVParam) }/{\partial \sysx}^{\Mytop}
)
(
\Drift(\sysx,\StoParam) 
+
\InMat(\sysx,\StoParam)
\InpOpt(\sysx, \InVecParam)
)
$, implying $\Bellman(\sysx,\preVParam,\InVecParam,\StoParam)=\BellA(\sysx)+\BellB(\sysx)$.
We consider the case of $\Input(\sysx(\MyT))=\InpOpt(\sysx(\MyT), \InVecParam)$.
For each $\StoParam$, there exists a locally unique solution $\sysx(\MyT,\StoParam)$ to the feedback system \eqref{eq:def_sys} because $\Drift(\sysx,\StoParam)$, $\InMat(\sysx,\StoParam)$, and $\InpOpt(\sysx,\InVecParam)$ are locally Lipschitz.
Because the global asymptotic stability implies the boundedness $\sup_{\MyT \geq 0 } \|\sysx(\MyT,\StoParam)\| < \infty$, there exists a global unique solution $\sysx(\MyT,\StoParam)$ on $\MyT\in [0,\infty)$ by using \cite[Proposition C.3.6]{Sontag98}.
Since $\InpOpt(\sysx,\InVecParam)$, ${\partial \EachValue(\sysx,\preVParam) }/{\partial \sysx}$, $\Drift(\sysx,\StoParam)$, $\InMat(\sysx,\StoParam)$, $\uCost(\sysx)$, and $\xCost(\sysx)$ are continuous in $\sysx$ and $\sysx(\MyT,\StoParam)$ is continuous in $\MyT$, $\BellA$ and $\BellB$ are integrable regarding $\MyT$.
In the following, the arguments $(\MyT,\StoParam) $ are often omitted for the brevity of the notation. 
For each $\StoParam$ and each $\sysx(0,\StoParam) \in\mathbb{R}^{\DimX}$, 
integrating $\Bellman(\sysx,\preVParam,\InVecParam,\StoParam)$ and using \eqref{eq:Bellman_residual_bound} yield
\begin{align}
	\int_{0}^{\terminalT}
	\big(
	\BellA(\sysx)+\BellB(\sysx)
	\big)
	\mathrm{d}\MyT
	&\leq
	\int_{0}^{\terminalT}
	\BellResiCoef
	\BellA(\sysx)
	\mathrm{d}\MyT
	.\label{eq:Bell_pf0}
\end{align}
For every $\terminalT>0$, this inequality indicates
\begin{align}
	\BellLim(\terminalT)
	:=
	\int_{0}^{\terminalT}
	\big(
	(\BellResiCoef-1)\BellA(\sysx)-\BellB(\sysx)
	\big)
	\mathrm{d}\MyT
	&
	\geq
	0
	.\label{eq:Bell_pf1}
\end{align}
Because of $\JBasis(0)=0$, the global asymptotic stability of the feedback system \eqref{eq:def_sys} leads to 
\begin{align}
	\lim_{\terminalT \to \infty} \EachValue(\sysx(\terminalT),\preVParam)
	&=
	\lim_{\sysx(\terminalT) \to 0}
	\JBasis(\sysx(\terminalT))^{\Mytop}  \preVParam
	=0
	.
\end{align}
For the given $\sysx(0,\StoParam) \in\mathbb{R}^{\DimX}$, this result indicates
\begin{align}
	\lim_{\terminalT \to \infty} 
	\int_{0}^{\terminalT}
	\BellB(\sysx)
	\mathrm{d}\MyT	
	&=
	\lim_{\terminalT \to \infty} 
	\EachValue(\sysx(\terminalT),\preVParam) - \EachValue(\sysx(0),\preVParam)
	\nonumber\\&
	=- \EachValue(\sysx(0),\preVParam)
	.\label{eq:Bell_pf2}
\end{align}
Meanwhile, using the definition \eqref{eq:def_cost} and the finiteness of $\costJ(\InpOpt(\NonArg, \InVecParam), \sysx(0),\StoParam)$ gives
\begin{align}
	\lim_{\terminalT \to \infty}
	\int_{0}^{\terminalT}
	\BellA(\sysx)
	\mathrm{d}\MyT
	=
	\costJ(\InpOpt(\NonArg, \InVecParam), \sysx(0),\StoParam)
	.\label{eq:Bell_pf3}
\end{align} 
Let us define
\begin{align}
	\BellExact
	:=
	(\BellResiCoef-1)
	\costJ(\InpOpt(\NonArg, \InVecParam), \sysx(0),\StoParam)
	+
	\EachValue(\sysx(0),\preVParam)
	.\label{eq:def_BellExact}	
\end{align}
Using \eqref{eq:Bell_pf1}, \eqref{eq:Bell_pf2}, and \eqref{eq:Bell_pf3} provides
\begin{align}
	&	
	\forall \PROOFepsilon>0,
	\exists \terminalT^{\prime}>0,
	\terminalT>\terminalT^{\prime}
	\Rightarrow
	|
	\BellLim(\terminalT)
	-
	\BellExact
	|
	<\PROOFepsilon
	.
\end{align}
Supposing that $\BellExact<0$ holds, setting $0< \PROOFepsilon < -\BellExact$ indicates 
\begin{align}
	\exists \terminalT^{\prime}>0,
	\terminalT>\terminalT^{\prime}
	\Rightarrow
	\BellLim(\terminalT)
	<
	\PROOFepsilon+\BellExact
	<0
	.
\end{align}
Because this contradicts \eqref{eq:Bell_pf1}, we have $\BellExact\geq 0$.
Combining this result with \eqref{eq:def_BellExact} yields
\begin{align}
	\costJ(\InpOpt(\NonArg, \InVecParam), \sysx(0),\StoParam)
	-
	\EachValue(\sysx(0),\preVParam)
	\leq \BellResiCoef 	\costJ(\InpOpt(\NonArg, \InVecParam), \sysx(0),\StoParam)
	.\label{eq:Bell_pf_result1}
\end{align}
In the same manner, we obtain
\begin{align}
	\costJ(\InpOpt(\NonArg, \InVecParam), \sysx(0),\StoParam)
	-
	\EachValue(\sysx(0),\preVParam)
	\geq -\BellResiCoef 	\costJ(\InpOpt(\NonArg, \InVecParam), \sysx(0),\StoParam)
	,\label{eq:Bell_pf_result2}
\end{align}
by starting from the following inequality instead of \eqref{eq:Bell_pf0}:
\begin{align}
	\int_{0}^{\terminalT}
	\big(
	\BellA(\sysx)+\BellB(\sysx)
	\big)
	\mathrm{d}\MyT
	&\geq
	-
	\int_{0}^{\terminalT}
	\BellResiCoef
	\BellA(\sysx)
	\mathrm{d}\MyT
	.
\end{align}
Using \eqref{eq:Bell_pf_result1} and \eqref{eq:Bell_pf_result2} yields \eqref{eq:Bellman_eq_equivalent_err} for each $\sysx(0,\StoParam) \in\mathbb{R}^{\DimX}$.
This completes the proof.

\section{Proof of Theorem \ref{thm:explicit_cost}}\label{pf:explicit_cost}

Substituting the parameterizations in \eqref{eq:def_parametric_input} and \eqref{eq:def_parametric_cost} and Assumption \ref{ass:polypoly_sys} into the Bellman residual \eqref{eq:Bellman_eq} yields
\begin{align}
	&\Bellman(\sysx,\preVParam,\InVecParam,\StoParam)
	\nonumber\\&
	=
	\xCost(\sysx)
	+
	\frac{1}{2}
	\InVecParam^{\Mytop} 
	\uBasis(\sysx)
	\uCost(\sysx)
	\uBasis(\sysx)^{\Mytop} 
	\InVecParam
	\nonumber\\&\quad
	+
	\preVParam^{\Mytop}  \frac{\partial \JBasis(\sysx)}{\partial \sysx^{\Mytop}}
	\sum_{\IDAffine=1}^{\NumAffine}  \sDrift{\IDAffine}(\StoParam) 
	\Big(
	\xDrift{\IDAffine}(\sysx) 
	+
	\xInMat{\IDAffine}(\sysx)
	\uBasis(\sysx)^{\Mytop} 
	\InVecParam
	\Big)
	. \label{eq:Bellman_eq_trans}
\end{align}
Using the properties of 
$\VEC{\NotationMat_{1}\NotationMat_{2}\NotationMat_{3}}=(\NotationMat_{3}^{\Mytop} \otimes \NotationMat_{1}) \VEC{\NotationMat_{2}}$ \citep[Eq. (3.106) in Sec. 3.2.10.2]{Gentle17} 
and
$(\NotationMat_{1}\otimes\NotationMat_{2})( \NotationMat_{3}\otimes\NotationMat_{4})=\NotationMat_{1}\NotationMat_{3}\otimes \NotationMat_{2}\NotationMat_{4}$ \citep[Eq. (3.101) in Sec. 3.2.10.2]{Gentle17}
for matrices $\NotationMat_{1}$, $\NotationMat_{2}$, $\NotationMat_{3}$, and $\NotationMat_{4}$,
we derive the relations:
\begin{align}
	&\xCost(\sysx)
	+
	\frac{1}{2}
	\InVecParam^{\Mytop} 
	\uBasis(\sysx)
	\uCost(\sysx)
	\uBasis(\sysx)^{\Mytop} 
	\InVecParam
	\nonumber\\&\;=
	\xCost(\sysx)
	+
	( \InVecParam^{\Mytop} \otimes \InVecParam^{\Mytop})
	{\VEC{  
			\uBasis(\sysx)	\uCost(\sysx)	\uBasis(\sysx)^{\Mytop}  /2
	}}	
	\nonumber\\&\;
	=
	\xBterm^{\Mytop}
	\begin{bmatrix}
		1 \\ \InVecParam \otimes \InVecParam
	\end{bmatrix}	
	,
	\\
	&
	\preVParam^{\Mytop}  \frac{\partial \JBasis(\sysx)}{\partial \sysx^{\Mytop}}
	\sum_{\IDAffine=1}^{\NumAffine}  \sDrift{\IDAffine}(\StoParam) 
	\Big(
	\xDrift{\IDAffine}(\sysx) 
	+
	\xInMat{\IDAffine}(\sysx)
	\uBasis(\sysx)^{\Mytop} 
	\InVecParam
	\Big)
	\nonumber\\&
	=
	\sum_{\IDAffine=1}^{\NumAffine}  \sDrift{\IDAffine}(\StoParam) 
	\preVParam^{\Mytop}  \frac{\partial \JBasis(\sysx)}{\partial \sysx^{\Mytop}}
	[ \xDrift{\IDAffine}(\sysx)  ,  \xInMat{\IDAffine}(\sysx) \uBasis(\sysx)^{\Mytop}   ]
	\begin{bmatrix}
		1 \\ \InVecParam
	\end{bmatrix}
	\nonumber\\&
	=
	\sum_{\IDAffine=1}^{\NumAffine}  \sDrift{\IDAffine}(\StoParam) 
	\Big(
	\begin{bmatrix}
		1 \\ \InVecParam
	\end{bmatrix}^{\Mytop}
	\otimes
	\preVParam^{\Mytop} 
	\Big)
	{\xAterm{\IDAffine}}
	\nonumber\\&
	=
	\sum_{\IDAffine=1}^{\NumAffine}  \sDrift{\IDAffine}(\StoParam) 
	{\xAterm{\IDAffine}}^{\Mytop}
	\Big(
	\begin{bmatrix}
		1 \\ \InVecParam
	\end{bmatrix}
	\otimes
	\preVParam
	\Big)
	\nonumber\\&
	=
	\sum_{\IDAffine=1}^{\NumAffine}  \sDrift{\IDAffine}(\StoParam) 
	{\xAterm{\IDAffine}}^{\Mytop}
	\Big( 
	\begin{bmatrix}
		1 \\ \InVecParam
	\end{bmatrix}
	\otimes 
	\Identity{\DimVParam} 
	\Big)
	( 1 \otimes  \preVParam )
	\nonumber\\&
	=
	\sum_{\IDAffine=1}^{\NumAffine}  \sDrift{\IDAffine}(\StoParam) 
	{\xAterm{\IDAffine}}^{\Mytop}
	\Big( 
	\begin{bmatrix}
		1 \\ \InVecParam
	\end{bmatrix}
	\otimes 
	\Identity{\DimVParam} 
	\Big)
	\preVParam
	.
\end{align}
Substituting these two relations into \eqref{eq:Bellman_eq_trans} yields
\begin{align}
	&\Bellman(\sysx,\preVParam,\InVecParam,\StoParam)
	\nonumber\\&
	=
	\xBterm^{\Mytop}
	\begin{bmatrix}
		1 \\ \InVecParam \otimes \InVecParam
	\end{bmatrix}	
	+
	\sum_{\IDAffine=1}^{\NumAffine}  \sDrift{\IDAffine}(\StoParam) 
	{\xAterm{\IDAffine}}^{\Mytop}
	\Big( 
	\begin{bmatrix}
		1 \\ \InVecParam
	\end{bmatrix}
	\otimes 
	\Identity{\DimVParam} 
	\Big)
	\preVParam
	.\label{eq:Bellman_eq_trans2}
\end{align}
Then, we obtain 
\begin{align}
	&\Bellman( \sysx, \preVParam,\InVecParam, \StoParam)^{2}
	\nonumber\\&
	=
	\preVParam^{\Mytop}
	\Bigg(
	\sum_{\IDAffine=1}^{\NumAffine}
	\sum_{\IDbAffine=1}^{\NumAffine}
	\sDrift{\IDAffine}(\StoParam) 
	\sDrift{\IDbAffine}(\StoParam)
	\Big( 
	\begin{bmatrix}
		1 \\ \InVecParam
	\end{bmatrix}
	\otimes 
	\Identity{\DimVParam} 
	\Big)^{\Mytop} 
	\nonumber\\&\quad\times
	{\xAterm{\IDAffine}}
	{\xAterm{\IDbAffine}}^{\Mytop}
	\Big( 
	\begin{bmatrix}
		1 \\ \InVecParam
	\end{bmatrix}
	\otimes 
	\Identity{\DimVParam} 
	\Big)
	\Bigg)
	\preVParam
	\nonumber\\&\quad
	+
	2 
	\preVParam^{\Mytop}
	\sum_{\IDAffine=1}^{\NumAffine}  \sDrift{\IDAffine}(\StoParam) 
	\Big( 
	\begin{bmatrix}
		1 \\ \InVecParam
	\end{bmatrix}
	\otimes 
	\Identity{\DimVParam} 
	\Big)^{\Mytop} 
	{\xAterm{\IDAffine}}
	\xBterm^{\Mytop}
	\begin{bmatrix}
		1 \\ \InVecParam \otimes \InVecParam
	\end{bmatrix}
	\nonumber\\&\quad
	+
	\Big(
	\xBterm^{\Mytop}
	\begin{bmatrix}
		1 \\ \InVecParam \otimes \InVecParam
	\end{bmatrix}
	\Big)^{2}
	.
\end{align}
Because all the functions of $\sysx$ are separated from the other components in $\Bellman( \sysx, \preVParam,\InVecParam, \StoParam)^{2}$, the objective function in \eqref{eq:def_OptVParam} is given as follows:
\begin{align}
	&\MarginalMap[\Big]{   \Bellman( \NonArg, \preVParam,\InVecParam, \StoParam)^{2}  }
	+
	{\CoefReguVParam} \|\preVParam\|^{2}
	\nonumber\\&
	=
	\preVParam^{\Mytop}
	\Bigg(
	\sum_{\IDAffine=1}^{\NumAffine}
	\sum_{\IDbAffine=1}^{\NumAffine}
	\sDrift{\IDAffine}(\StoParam) 
	\sDrift{\IDbAffine}(\StoParam)
	{\newDDmat{\InVecParam}{\IDAffine}{\IDbAffine}}
	\Bigg)
	\preVParam
	\nonumber\\&\quad
	+
	2 
	\preVParam^{\Mytop}
	\sum_{\IDAffine=1}^{\NumAffine}  \sDrift{\IDAffine}(\StoParam) 
	{\newDNmat{\InVecParam}{\IDAffine}} 
	+
	\MarginalMap[\Big]{
		\Big(
		\xBterm^{\Mytop}
		\begin{bmatrix}
			1 \\ \InVecParam \otimes \InVecParam
		\end{bmatrix}
		\Big)^{2}
	}
	, \label{eq:pf_Bellman}
\end{align}
because
${\CoefReguVParam}
\Identity{\DimVParam} 
=
\sum_{\IDAffine=1}^{\NumAffine}
\sum_{\IDbAffine=1}^{\NumAffine}
\sDrift{\IDAffine}(\StoParam) 
\sDrift{\IDbAffine}(\StoParam) 
{\CoefReguVParam}
\Identity{\DimVParam} 
$ holds by \eqref{eq:sum_sDrift_one}.
All the terms in \eqref{eq:pf_Bellman} are bounded and thus well-defined because 
using Cauchy-Schwarz inequality provides
\begin{align}
\MarginalMap[]{ 
	{\El{\xAterm{\IDAffine}}{\IDEl}}
	{\El{\xAterm{\IDbAffine}}{\IDbEl}}
}^{2}
&\leq
\MarginalMap[]{ 
	{\El{\xAterm{\IDAffine}}{\IDEl}}^{2}
}
\MarginalMap[]{ 
	{\El{\xAterm{\IDbAffine}}{\IDbEl}}^{2}
}
\nonumber\\&
\leq
\MarginalMap[]{ 
	\| {\xAterm{\IDAffine}}\|^{2} 
}
\MarginalMap[]{ 
	\| {\xAterm{\IDbAffine}}\|^{2} 
}
\nonumber\\&
< \infty
.
\end{align}
In a manner of standard quadratic minimization, the optimal parameter $\VParam(\InVecParam,\StoParam)$ in \eqref{eq:def_OptVParam} is given by \eqref{eq:solution_OptVParam}.
This completes the proof.

\section{Proof of Theorem \ref{thm:explicit_grad}}\label{pf:explicit_grad}

We  derive the partial derivative of $\Expect{\IniState}[ \EachValue(\IniState,\VParam(\InVecParam,\StoParam))] $.
For the simplicity of the notation, let us denote
\begin{align}
	\StoDDmat(\InVecParam,\StoParam)
	&:=
	\sum_{\IDAffine=1}^{\NumAffine}
	\sum_{\IDbAffine=1}^{\NumAffine}
	\sDrift{\IDAffine}(\StoParam) 
	\sDrift{\IDbAffine}(\StoParam)
	{\newDDmat{\InVecParam}{\IDAffine}{\IDbAffine}}
	,
	\\
	\StoDNmat(\InVecParam,\StoParam)
	&:= 
	\sum_{\IDAffine=1}^{\NumAffine}
	\sDrift{\IDAffine}(\StoParam)
	{\newDNmat{\InVecParam}{\IDAffine}} 
	,
\end{align}
which implies $\VParam(\InVecParam,\StoParam)= - \StoDDmat(\InVecParam,\StoParam)^{-1}\StoDNmat(\InVecParam,\StoParam)$.
Using the property of
$\partial \NotationMat^{-1}/ \partial {\El{\InVecParam}{\IDbEl}}
=
- \NotationMat^{-1}
(\partial \NotationMat/ \partial {\El{\InVecParam}{\IDbEl}})
\NotationMat^{-1}
$ for a matrix-valued function $\NotationMat$ \citep[Eq. (C.21)]{Bishop06},
the partial derivative of
$\Expect{\IniState}[ \EachValue(\IniState,\VParam(\InVecParam,\StoParam))] = \Expect{\IniState}[   
\JBasis(\IniState) 
]^{\Mytop} 
\VParam(\InVecParam,\StoParam) $
is obtained
as follows:
\begin{align}
	&
	\frac{\partial }{\partial \El\InVecParam{\IDbEl}}
	\Expect{\IniState}[    
	\EachValue(\IniState,\VParam(\InVecParam,\StoParam))
	]   		
	\nonumber\\&
	=
	\Expect{\IniState}[   
	\JBasis(\IniState) 
	]^{\Mytop} 
	\frac{\partial }{\partial \El\InVecParam{\IDbEl}}
	\VParam(\InVecParam,\StoParam)
	\nonumber\\&
	=
	\Expect{\IniState}[   
	\JBasis(\IniState) 
	]^{\Mytop} 
	\Big(
	-
	\StoDDmat(\InVecParam,\StoParam)^{-1} 
	\frac{\partial \StoDNmat(\InVecParam,\StoParam) }{\partial \El\InVecParam{\IDbEl}}
	\nonumber\\&\quad
	+
	\StoDDmat(\InVecParam,\StoParam)^{-1} 
	\frac{\partial \StoDDmat(\InVecParam,\StoParam)}{\partial \El\InVecParam{\IDbEl}}
	\StoDDmat(\InVecParam,\StoParam)^{-1} 
	\StoDNmat(\InVecParam,\StoParam)
	\Big)
	\nonumber\\&
	=
	\Expect{\IniState}[   
	\JBasis(\IniState) 
	]^{\Mytop} 
	\Big(
	-
	\StoDDmat(\InVecParam,\StoParam)^{-1} 
	\frac{\partial \StoDNmat(\InVecParam,\StoParam) }{\partial \El\InVecParam{\IDbEl}}
	\nonumber\\&\quad
	-
	\StoDDmat(\InVecParam,\StoParam)^{-1} 
	\frac{\partial \StoDDmat(\InVecParam,\StoParam)}{\partial \El\InVecParam{\IDbEl}}
	\VParam(\InVecParam,\StoParam)
	\Big)
	\nonumber\\&
	=
	-
	\Expect{\IniState}[   
	\JBasis(\IniState) 
	]^{\Mytop} 
	\StoDDmat(\InVecParam,\StoParam)^{-1} 
	\Big(		
	\frac{\partial \StoDNmat(\InVecParam,\StoParam) }{\partial \El\InVecParam{\IDbEl}}
	+
	\frac{\partial \StoDDmat(\InVecParam,\StoParam)}{\partial \El\InVecParam{\IDbEl}}
	\VParam(\InVecParam,\StoParam)
	\Big)
	. 
\end{align}
This result corresponds to \eqref{eq:explicit_cost_derivative}.
This completes the proof.

\section{Proof of Lemma \ref{thm:design_penaF}}\label{pf:design_penaF}

As \eqref{eq:def_symmetricSOSdotVMat} implies that ${\SOSdotVMat{\IDAffine}}$ is symmetric,
the half vectorization ${\VECH{  
		{\SOSdotVMat{\IDAffine}}
}}$ is well-defined.
Using the property of 
$\VEC{\NotationMat_{1}\NotationMat_{2}\NotationMat_{3}}=(\NotationMat_{3}^{\Mytop} \otimes \NotationMat_{1}) \VEC{\NotationMat_{2}}$ \citep[Eq. (3.106) in Sec. 3.2.10.2]{Gentle17}, we have
\begin{align}
	\SOSdotVBasis(\sysx)^{\Mytop}
	{\SOSdotVMat{\IDAffine}}
	\SOSdotVBasis(\sysx)
	&=
	(\SOSdotVBasis(\sysx)^{\Mytop} \otimes \SOSdotVBasis(\sysx)^{\Mytop} )
	{\VEC{  
			{\SOSdotVMat{\IDAffine}} 
	}}
	\nonumber\\&
	=
	{\VEC{  \SOSdotVBasis(\sysx) \SOSdotVBasis(\sysx)^{\Mytop} }}^{\Mytop}
	{\VEC{  
			{\SOSdotVMat{\IDAffine}}
	}}
	\nonumber\\&
	=
	{\VECH{  \SOSdotVBasis(\sysx) \SOSdotVBasis(\sysx)^{\Mytop} }}^{\Mytop}
	\connectMatA
	{\VECH{  
			{\SOSdotVMat{\IDAffine}}
	}}
	,\label{eq:SOSdotVMat_pf1}
\end{align}
where $\connectMatA$ is a positive definite diagonal matrix with diagonal components equal to $1$ or $2$.

Let $\SOSdotVBBasis(\sysx) \in \mathbb{R}^{\DimSOSdotVBBasis}$ be 
the non-redundant form of $ {\VECH{  \SOSdotVBasis(\sysx) \SOSdotVBasis(\sysx)^{\Mytop} }}$,
where 
$\DimSOSdotVBBasis\leq\DimSOSdotVBasis(\DimSOSdotVBasis+1)/2 $ holds clearly.
Then, there exists a matrix $\connectMatB \in \{0,1\}^{\DimSOSdotVBBasis\times \DimSOSdotVBasis(\DimSOSdotVBasis+1)/2  }$ 
satisfying
\begin{align}
	{\VECH{  \SOSdotVBasis(\sysx) \SOSdotVBasis(\sysx)^{\Mytop} }} 
	=
	\connectMatB^{\Mytop} 
	\SOSdotVBBasis(\sysx)
	.\label{eq:def_connectMatB}
\end{align}
Substituting this relation into \eqref{eq:SOSdotVMat_pf1} yields
\begin{align}
	\SOSdotVBasis(\sysx)^{\Mytop}
	{\SOSdotVMat{\IDAffine}}
	\SOSdotVBasis(\sysx)
	&=
	\SOSdotVBBasis(\sysx)^{\Mytop}
	\connectMatB
	\connectMatA
	{\VECH{  
			{\SOSdotVMat{\IDAffine}}
	}}
	.\label{eq:SOSdotVMat_pf2}
\end{align}
Next, 
using Assumption \ref{ass:SOS} {\MyLabeldotVpoly} yields
\begin{align}
&
\SOSBasis(\sysx)^{\Mytop}
{\SOSdotVpreMat{\IDAffine}}(\sysx,\LyapMat,\InVecParam)
\SOSBasis(\sysx)
\nonumber\\&
=
\sum_{\IDEl=1}^{\DimSOSBasis}
\sum_{\IDbEl=1}^{\DimSOSBasis}
{\El{ {\SOSdotVpreMat{\IDAffine}}(\sysx,\LyapMat,\InVecParam) }{\IDEl,\IDbEl}}
{\El{ \SOSBasis(\sysx) }{\IDEl}}
{\El{ \SOSBasis(\sysx) }{\IDbEl}}	
\nonumber\\&	
=
\sum_{\IDEl=1}^{\DimSOSBasis}
\sum_{\IDbEl=1}^{\DimSOSBasis}
{\coeffdotVpreSqrtMat{\IDEl}{\IDbEl}}^{\Mytop} 
( {\El{ \SOSBasis(\sysx) }{\IDEl}} \SOSmatBasis(\sysx)
\otimes 
  {\El{ \SOSBasis(\sysx) }{\IDbEl}}	\SOSmatBasis(\sysx) ) 
.\label{eq:SOSdotVMat_pf3_1}
\end{align}
Because the components of $\SOSdotVBasis(\sysx)$ contain ${\El{ \SOSBasis(\sysx) }{\IDEl}} \SOSmatBasis(\sysx)$ and ${\El{ \SOSBasis(\sysx) }{\IDbEl}} \SOSmatBasis(\sysx)$,
there exists $\coeffVECdotVpreSqrtMat  \in \mathbb{R}^{\DimSOSdotVBBasis}$ satisfying 
\begin{align}
	\SOSBasis(\sysx)^{\Mytop}
	{\SOSdotVpreMat{\IDAffine}}(\sysx,\LyapMat,\InVecParam)
	\SOSBasis(\sysx)
	=
	\SOSdotVBBasis(\sysx)^{\Mytop}  \coeffVECdotVpreSqrtMat
	.\label{eq:SOSdotVMat_pf3}
\end{align}
where $\coeffVECdotVpreSqrtMat$ depends on $(\InVecParam,\LyapMat)$ and is uniquely determined because $\SOSdotVBBasis(\sysx)$ is the non-redundant form.
Substituting \eqref{eq:SOSdotVMat_pf2} and \eqref{eq:SOSdotVMat_pf3} into \eqref{eq:def_dotLyapF_new}	yields
\begin{align}
	\SOSdotVBBasis(\sysx)^{\Mytop}
	(
	\connectMatB \connectMatA
	{\VECH{  
			{\SOSdotVMat{\IDAffine}}
	}}
	-
	\coeffVECdotVpreSqrtMat
	)
	=0
	,\;
	\forall \sysx 
	.\label{eq:def_dotLyapF_anotherA}
\end{align}
This zero function must be a zero polynomial \citep[Proposition 5]{Cox07}.
Thus,  \eqref{eq:def_dotLyapF_anotherA} and \eqref{eq:def_dotLyapF_new} are equivalent to the following linear equation:
\begin{align}
	\connectMatB \connectMatA
	{\VECH{  
			{\SOSdotVMat{\IDAffine}}
	}}
	=
	\coeffVECdotVpreSqrtMat
	.\label{eq:def_dotLyapF_anotherB}
\end{align}

Here, let us consider the case of $\DimSOSdotVBBasis=\DimSOSdotVBasis(\DimSOSdotVBasis+1)/2 $, that is, $\connectMatB$ is a square permutation matrix that is  orthogonal and thus nonsingular.
Because both $\connectMatA$ and $\connectMatB$ are nonsingular,
$
{\VECH{  
		{\SOSdotVMat{\IDAffine}}
}}
=  \connectMatA^{-1} \connectMatB^{-1}  \coeffVECdotVpreSqrtMat
$ is the unique solution to \eqref{eq:def_dotLyapF_anotherB}.
This is the unique solution satisfying \eqref{eq:def_dotLyapF_new}--\eqref{eq:def_symmetricSOSdotVMat} under the settings of $\AScoefmat=0$ and $\AScoefscalar=0$.

Meanwhile, we consider the other case: $\DimSOSdotVBBasis<\DimSOSdotVBasis(\DimSOSdotVBasis+1)/2 $.
Because of \eqref{eq:def_connectMatB}, we have
$\sum_{\IDEl=1}^{\DimSOSdotVBBasis} {\El{  \connectMatB  }{\IDEl,\IDbEl}}=1$ for each $\IDbEl$
and
$\sum_{\IDbEl=1}^{\DimSOSdotVBasis(\DimSOSdotVBasis+1)/2} {\El{  \connectMatB  }{\IDEl,\IDbEl}}\geq 1$ for each $\IDEl$.
Thus, there exists a permutation matrix $\connectMatC\in \{0,1\}^{\DimSOSdotVBasis(\DimSOSdotVBasis+1)/2 \times \DimSOSdotVBasis(\DimSOSdotVBasis+1)/2}$ and
a matrix $ \connectMatD \in \{0,1\}^{ \DimSOSdotVBBasis \times (\DimSOSdotVBasis(\DimSOSdotVBasis+1)/2 - \DimSOSdotVBBasis) }$ such that
we have $\connectMatB = [{\Identity{\DimSOSdotVBBasis}} , \connectMatD] \connectMatC$.
Substituting this into \eqref{eq:def_dotLyapF_anotherB} provides
\begin{align}
	[{\Identity{\DimSOSdotVBBasis}} , \connectMatD] \connectMatC
	\connectMatA
	{\VECH{  
			{\SOSdotVMat{\IDAffine}}  
	}}
	=
	\coeffVECdotVpreSqrtMat
	.\label{eq:def_dotLyapF_anotherC}
\end{align}
By setting $\AScoefscalar=1$ and 
$
\AScoefmat 
=
[ 0 , {\Identity{(\DimSOSdotVBasis(\DimSOSdotVBasis+1)/2 - \DimSOSdotVBBasis)}} ]  \connectMatC \connectMatA
$, 
combining \eqref{eq:def_dotLyapF_anotherC} and \eqref{eq:def_ASconstraint} yields
\begin{align}
	\connectMatE
	\connectMatC
	\connectMatA
	{\VECH{  
			{\SOSdotVMat{\IDAffine}}
	}}
	&=
	\begin{bmatrix}
		\coeffVECdotVpreSqrtMat
		\\
		{\ASeach{\IDAffine}}
	\end{bmatrix}
	,\label{eq:def_dotLyapF_anotherD}
	\\
	\connectMatE
	&:=
	\begin{bmatrix}
		{\Identity{\DimSOSdotVBBasis}} & \connectMatD\\
		0 &  {\Identity{(\DimSOSdotVBasis(\DimSOSdotVBasis+1)/2 - \DimSOSdotVBBasis)}} 
	\end{bmatrix}
	.
\end{align}
The matrix $\connectMatE$ is nonsingular because it is a triangular matrix with non-zero diagonal components.
Because $\connectMatE$, $\connectMatC$, and $\connectMatA$ are nonsingular, 
$
{\VECH{  
		{\SOSdotVMat{\IDAffine}}
}}
=  \connectMatA^{-1} \connectMatC^{-1}   \connectMatE^{-1} 
[
\coeffVECdotVpreSqrtMat^{\Mytop}
,
{\ASeach{\IDAffine}}^{\Mytop}
]^{\Mytop}
$ is the unique solution to \eqref{eq:def_dotLyapF_anotherD}.
This is the unique solution satisfying \eqref{eq:def_dotLyapF_new}--\eqref{eq:def_symmetricSOSdotVMat}.
Note that $\connectMatA$ and $\connectMatC$ are determined from the monomial bases, implying that $\AScoefmat$ is independent of $\IDAffine$.

Finally, $\coeffVECdotVpreSqrtMat$ is bilinear in $(\InVecParam,\LyapMat)$ because of the definition \eqref{eq:def_SOSdotVpreMat} of ${\SOSdotVpreMat{\IDAffine}}(\sysx,\LyapMat,\InVecParam)$. 
Thus,  the obtained solution to \eqref{eq:def_dotLyapF_anotherB} or \eqref{eq:def_dotLyapF_anotherD} is bilinear in $(\InVecParam,\LyapMat)$ and linear in ${\ASeach{\IDAffine}}$.
This completes the proof.

\section{Proof of Theorem \ref{thm:stabilizing_optim}}\label{pf:stabilizing_optim}

Firstly, we show that the following conditions hold for a given $\IDIteGrad> 0$:
\begin{align}
	{\Ite{\LyapMat}{\IDIteGrad}}+{\Ite{\LyapMat}{\IDIteGrad}}^{\Mytop} \succ 0
	, \label{eq:cond_robust_stab_P_all}
	\\
	{\SOSdotVMat{\IDAffine}}({\Ite{\InVecParam}{\IDIteGrad}},{\Ite{\LyapMat}{\IDIteGrad}},{\Ite{\ASparam}{\IDIteGrad}}) \succ 0
	. \label{eq:cond_robust_stab_T_all}
\end{align}
Supposing that \eqref{eq:cond_robust_stab_P_all} or \eqref{eq:cond_robust_stab_T_all} is not satisfied,
\eqref{eq:def_PenaltyFuncRegu} indicates $\PenaltyFuncRegu({\Ite{\InVecParam}{\IDIteGrad}},{\Ite{\LyapMat}{\IDIteGrad}},{\Ite{\ASparam}{\IDIteGrad}})=\iniPenaltyFuncRegu$.
Meanwhile, using \eqref{eq:cond_iniPenaltyFuncRegu}, \eqref{pf:ObjF_mono_dec}, and \eqref{eq:cost_morethan_LB} yields 
$
\appcostLB
+
\PenaltyFuncRegu({\Ite{\InVecParam}{\IDIteGrad}},{\Ite{\LyapMat}{\IDIteGrad}},{\Ite{\ASparam}{\IDIteGrad}})
\leq
\ObjF(  {\Ite{\InVecParam}{\IDIteGrad}}   ,{\Ite{\LyapMat}{\IDIteGrad}} ,{\Ite{\ASparam}{\IDIteGrad}})
\leq
\ObjF(  {\Ite{\InVecParam}{0}}   ,{\Ite{\LyapMat}{0}} ,{\Ite{\ASparam}{0}})
<
\appcostLB
+
\iniPenaltyFuncRegu
$, 
indicating $\PenaltyFuncRegu({\Ite{\InVecParam}{\IDIteGrad}},{\Ite{\LyapMat}{\IDIteGrad}},{\Ite{\ASparam}{\IDIteGrad}})<\iniPenaltyFuncRegu$.
Because this contradicts  $\PenaltyFuncRegu({\Ite{\InVecParam}{\IDIteGrad}},{\Ite{\LyapMat}{\IDIteGrad}},{\Ite{\ASparam}{\IDIteGrad}})=\iniPenaltyFuncRegu$, \eqref{eq:cond_robust_stab_P_all} and \eqref{eq:cond_robust_stab_T_all} hold  for the given $\IDIteGrad\geq 0$.
Here, the assumptions \eqref{eq:cond_robust_stab_P} and  \eqref{eq:cond_robust_stab_T} are necessary for the case of $\IDIteGrad=0$.
If they are not assumed,
\eqref{eq:cond_iniPenaltyFuncRegu} indicates $\Expect{\StoParam}[\Expect{\IniState}[    
\EachValue(\IniState,\VParam(  {\Ite{\InVecParam}{0}}    ,\StoParam))
] ]+ \iniPenaltyFuncRegu < \appcostLB
+
\iniPenaltyFuncRegu $ that contradicts \eqref{eq:cost_morethan_LB} with $\IDIteGrad=0$.

Next, we prove the robust global asymptotic stability by using \eqref{eq:cond_robust_stab_P_all} and \eqref{eq:cond_robust_stab_T_all}.
We set
$(\InVecParam, \LyapMat,\ASparam)=({\Ite{\InVecParam}{\IDIteGrad}},{\Ite{\LyapMat}{\IDIteGrad}},{\Ite{\ASparam}{\IDIteGrad}})$, recalling the Lyapunov function $\LyapF(\sysx):= \SOSBasis(\sysx)^{\Mytop} \LyapMat\SOSBasis(\sysx) 
= \SOSBasis(\sysx)^{\Mytop}( \LyapMat + \LyapMat^{\Mytop} )\SOSBasis(\sysx) /2 $ in \eqref{eq:def_Lyap}.
Then, $\LyapF(\sysx)$ is positive definite because  ${\Ite{\LyapMat}{\IDIteGrad}}+{\Ite{\LyapMat}{\IDIteGrad}}^{\Mytop} \succ 0$ holds by \eqref{eq:cond_robust_stab_P_all} and $\SOSBasis(\sysx)$ is strict, that is, $\SOSBasis(\sysx)= 0 \Leftrightarrow \sysx= 0$.

Here,
the strict $\SOSBasis(\sysx)$ contains the components $\El{\sysx}{\IDEl}^{\El{\degPow}{\IDEl}}$ for all $\IDEl$ and for some $\degPow \in \{0,1,2,\dots\}^{\DimX}$ because of $\SOSBasis(\sysx)\neq 0$ even for $\sysx=[0,\dots,0,1,0,\dots,0]^{\Mytop}$.
This implies that we have $\|\SOSBasis(\sysx)\|\to \infty$ as $\|\sysx\|\to \infty$.
Thus,  $\LyapF(\sysx) \geq {\EigMin{   \LyapMat + \LyapMat^{\Mytop}    }} \|\SOSBasis(\sysx)\|^{2}/2 $ is radially unbounded.

In addition, using the property of 
$\VEC{\NotationMat_{1}\NotationMat_{2}\NotationMat_{3}}=(\NotationMat_{3}^{\Mytop} \otimes \NotationMat_{1}) \VEC{\NotationMat_{2}}$ \citep[Eq. (3.106) in Sec. 3.2.10.2]{Gentle17}, 
combining \eqref{eq:def_InBasis} with \eqref{eq:def_parametric_input} gives 
\begin{align}
\InpOpt(\sysx, \InVecParam) 
= 
\uBasis(\sysx)^{\Mytop} \InVecParam
&=
( \SOSBasis(\sysx)^{\Mytop}  \otimes \InBasis(\sysx)   )\InVecParam
\nonumber\\&=
\InBasis(\sysx)  {\InvVEC{\InVecParam}{\DimRowInBasis}{\DimSOSBasis}}   \SOSBasis(\sysx)
.  \label{eq:def_parametric_input_transform}
\end{align}
Substituting \eqref{eq:apply_FBcontrol}, \eqref{eq:def_PolyDriftMat}, and \eqref{eq:def_parametric_input_transform} into \eqref{eq:def_sys} yields
\begin{align}
\dotState
&
= \Drift(\sysx,\StoParam) 
+ \InMat(\sysx,\StoParam) 
\InpOpt(\sysx, \InVecParam)
\nonumber\\&
=
\sum_{\IDAffine=1}^{\NumAffine}
{\sDrift{\IDAffine}}(\StoParam)
\Big(
{\PolyDriftMat{\IDAffine}}(\sysx)
+ 
\InMat(\sysx,\StoParam) 
\InBasis(\sysx)  {\InvVEC{\InVecParam}{\DimRowInBasis}{\DimSOSBasis}} 
\Big)
\SOSBasis(\sysx)
. \label{eq:sys_poly_sos_representation}
\end{align}
Using this relation, \eqref{eq:def_SOSdotVpreMat}, and \eqref{eq:def_dotLyapF_new}, we obtain
\begin{align}
\dot{\LyapF}(\sysx,\StoParam) 
&=
\SOSBasis(\sysx)^{\Mytop} \LyapMat \frac{\partial \SOSBasis(\sysx)}{\partial \sysx^{\Mytop} } \dotState
+
\dotState^{\Mytop} \Big(\frac{\partial \SOSBasis(\sysx)}{\partial \sysx^{\Mytop} }\Big)^{\Mytop} \LyapMat \SOSBasis(\sysx)
\nonumber\\&
=
- 2 \sum_{\IDAffine=1}^{\NumAffine}  \sDrift{\IDAffine}(\StoParam) 
\SOSBasis(\sysx)^{\Mytop} 
{\SOSdotVpreMat{\IDAffine}}(\sysx,\LyapMat,\InVecParam)
\SOSBasis(\sysx)
\nonumber\\&
=
-2\sum_{\IDAffine=1}^{\NumAffine}
{\sDrift{\IDAffine}}(\StoParam)
\SOSdotVBasis(\sysx)^{\Mytop}
{\SOSdotVMat{\IDAffine}}(\InVecParam, \LyapMat,\ASparam)
\SOSdotVBasis(\sysx)
\nonumber\\&
\leq
-2\sum_{\IDAffine=1}^{\NumAffine}
{\sDrift{\IDAffine}}(\StoParam)
\min_{\IDbAffine \in \{1,2,\dots,\NumAffine\}}
\SOSdotVBasis(\sysx)^{\Mytop}
{\SOSdotVMat{\IDbAffine}}(\InVecParam, \LyapMat,\ASparam)
\SOSdotVBasis(\sysx)
\nonumber\\&
=
-2
\min_{\IDbAffine \in \{1,2,\dots,\NumAffine\}}
\SOSdotVBasis(\sysx)^{\Mytop}
{\SOSdotVMat{\IDbAffine}}(\InVecParam, \LyapMat,\ASparam)
\SOSdotVBasis(\sysx)
. \label{eq:Lyap_calculationA}
\end{align}	
Thus,
$- \dot{\LyapF}(\sysx,\StoParam) $ is positive definite because \eqref{eq:cond_robust_stab_T_all} implies
${\SOSdotVMat{\IDAffine}}(\InVecParam, \LyapMat,\ASparam)\succ 0$ and $\SOSdotVBasis(\sysx)$ is strict, meaning $\SOSdotVBasis(\sysx)= 0 \Leftrightarrow \sysx= 0$.
Furthermore, this positive definiteness is global, that is, $\lim_{\|\sysx\| \to \infty} \dot{\LyapF}(\sysx,\StoParam) \neq 0 $, because $\|\SOSdotVBasis(\sysx)\|\to \infty$ as $\|\sysx\|\to \infty$.
By virtue of the radially unbounded positive definiteness of $\LyapF(\sysx) $ and 
the global positive definiteness of $- \dot{\LyapF}(\sysx,\StoParam) $ for every $\StoParam$,
the robust global asymptotic stability of the feedback system \eqref{eq:def_sys} with $\Input(\sysx)=\InpOpt(\sysx, \InVecParam)$ is guaranteed.
This completes the proof.

\section{Proof of Theorem \ref{thm:stabilizing_initial}}\label{pf:stabilizing_initial}

Firstly, we show the statement {\MyLabelFeasibleWP}.
From Assumption \ref{ass:SOS} {\MyLabelSOS},
there exist $\LyapMat \succ 0$, $\InVecParam \in   \mathbb{R}^{\DimInVecParam} $, and ${\SOSpdMat{\IDAffine}} \succ 0$ satisfying \eqref{eq:SOS_ass}.
Let us set 
\begin{align}
\LyapInv&=\LyapMat^{-1} \succ 0
, \label{eq:set_LyapInv}
\\
\tmpInpParam&={\InvVEC{\InVecParam}{\DimRowInBasis}{\DimSOSBasis}} \LyapMat^{-1}
. \label{eq:set_tmpInpParam}
\end{align}
Because of \eqref{eq:def_SOSdotVpreMat} and \eqref{eq:def_tdLyapMat}, we have
\begin{align}
{\tdLyapMat{\IDAffine}}(\sysx,\LyapInv, \tmpInpParam)
= \LyapMat^{-1} {\SOSdotVpreMat{\IDAffine}}(\sysx,\LyapMat,\InVecParam) \LyapMat^{-1}
= \LyapInv  {\SOSdotVpreMat{\IDAffine}}(\sysx,\LyapMat,\InVecParam)  \LyapInv
\end{align}	
Thus, using \eqref{eq:SOS_ass} yields	
\begin{align}
&
{\tdLyapMat{\IDAffine}}(\sysx,\LyapInv, \tmpInpParam)
+{\tdLyapMat{\IDAffine}}(\sysx,\LyapInv, \tmpInpParam)^{\Mytop}
\nonumber\\&
= 
\LyapInv
(\SOSmatBasis(\sysx) \otimes {\Identity{\DimSOSBasis}} )^{\Mytop}
{\SOSpdMat{\IDAffine}} 
(\SOSmatBasis(\sysx) \otimes {\Identity{\DimSOSBasis}} )
\LyapInv
.
\end{align}	
Applying the property $(\NotationMat_{1}\otimes\NotationMat_{2})( \NotationMat_{3}\otimes\NotationMat_{4})=\NotationMat_{1}\NotationMat_{3}\otimes \NotationMat_{2}\NotationMat_{4}$ \citep[Eq. (3.101) in Sec. 3.2.10.2]{Gentle17}
for matrices $\NotationMat_{1}$, $\NotationMat_{2}$, $\NotationMat_{3}$, and $\NotationMat_{4}$,
we obtain
\begin{align}
(\SOSmatBasis(\sysx) \otimes {\Identity{\DimSOSBasis}} )
\LyapInv
&=
(\SOSmatBasis(\sysx) \otimes {\Identity{\DimSOSBasis}} )
(1 \otimes \LyapInv	)
\nonumber\\&
=
(\SOSmatBasis(\sysx) \times 1 )
\otimes
({\Identity{\DimSOSBasis}} \LyapInv	)
\nonumber\\&
=
({\Identity{\DimSOSmatBasis}}   \SOSmatBasis(\sysx))
\otimes
(\LyapInv {\Identity{\DimSOSBasis}} 	)
\nonumber\\&
=
\blockLyapInv
( \SOSmatBasis(\sysx) \otimes {\Identity{\DimSOSBasis}} 	)
, \label{eq:kron_transform}
\end{align}	
where $\blockLyapInv:=\Identity{\DimSOSmatBasis} \otimes \LyapInv$.
Therefore,
\begin{align}
&
{\tdLyapMat{\IDAffine}}(\sysx,\LyapInv, \tmpInpParam)
+{\tdLyapMat{\IDAffine}}(\sysx,\LyapInv, \tmpInpParam)^{\Mytop}
\nonumber\\&
= 
(\SOSmatBasis(\sysx) \otimes   {\Identity{\DimSOSBasis}} )^{\Mytop}
\blockLyapInv^{\Mytop}
{\SOSpdMat{\IDAffine}} 
\blockLyapInv
(\SOSmatBasis(\sysx) \otimes   {\Identity{\DimSOSBasis}} )
\end{align}		
Because  $\LyapInv \succ 0$ and ${\SOSpdMat{\IDAffine}}  \succ 0$ hold, we have $\blockLyapInv \succ 0$ and thus
$\blockLyapInv^{\Mytop}
{\SOSpdMat{\IDAffine}} 
\blockLyapInv \succ 0$.
By setting  
\begin{align}
\forall \IDAffine \in \{1,2,\dots,\NumAffine\},\quad
{\SOSRtdLyapMat{\IDAffine}} 
=
\blockLyapInv^{\Mytop}
{\SOSpdMat{\IDAffine}} 
\blockLyapInv
\succ 0
, \label{eq:set_SOSRtdLyapMat}
\end{align}	
we satisfy the first constraint in the SDP \eqref{eq:SDP_to_LyapInv_tmpInpParam}.
In addition, the following setting can be used:
\begin{align}
\SOSAepsilon
&=
\min \{
{\EigMin{  \LyapInv   }},
{\EigMin{  {\SOSRtdLyapMat{1}}   }},
\dots,
{\EigMin{  {\SOSRtdLyapMat{\NumAffine}}   }}
\}
> 0
. \label{eq:set_SOSAepsilon}
\end{align}
Thus, \eqref{eq:set_LyapInv}, \eqref{eq:set_tmpInpParam}, \eqref{eq:set_SOSRtdLyapMat}, and \eqref{eq:set_SOSAepsilon} are feasible solutions to the SDP \eqref{eq:SDP_to_LyapInv_tmpInpParam}.

Next, we show the statement {\MyLabelFeasibleAddP}.
From \eqref{eq:def_SOSdotVpreMat} and \eqref{eq:def_tdLyapMat}, we have
\begin{align}
{\SOSdotVpreMat{\IDAffine}}(\sysx,{\Ite{\LyapMat}{0}},{\Ite{\InVecParam}{0}})
&= \solLyapInv^{-1} {\tdLyapMat{\IDAffine}}(\sysx,\solLyapInv, \soltmpInpParam) \solLyapInv^{-1}
\nonumber\\&
= {\Ite{\LyapMat}{0}}  {\tdLyapMat{\IDAffine}}(\sysx,\solLyapInv, \soltmpInpParam)   {\Ite{\LyapMat}{0}}
\end{align}
Substituting this relation and a feasible solution ${\solSOSRtdLyapMat{\IDAffine}}$ into the first constraint of the SDP \eqref{eq:SDP_to_LyapInv_tmpInpParam} yields
\begin{align}
&
{\SOSdotVpreMat{\IDAffine}}(\sysx,{\Ite{\LyapMat}{0}},{\Ite{\InVecParam}{0}})
+
{\SOSdotVpreMat{\IDAffine}}(\sysx,{\Ite{\LyapMat}{0}},{\Ite{\InVecParam}{0}})^{\Mytop}
\nonumber\\&
= 
{\Ite{\LyapMat}{0}} 
(\SOSmatBasis(\sysx) \otimes   {\Identity{\DimSOSBasis}} )^{\Mytop}
{\solSOSRtdLyapMat{\IDAffine}}
(\SOSmatBasis(\sysx) \otimes   {\Identity{\DimSOSBasis}} )
{\Ite{\LyapMat}{0}}
.\label{eq:SOSdotVpreMat_pd1}
\end{align}
In the same manner as \eqref{eq:kron_transform}, using $\blockLyapMat:=\Identity{\DimSOSmatBasis} \otimes {\Ite{\LyapMat}{0}}$ provides
\begin{align}
(\SOSmatBasis(\sysx) \otimes {\Identity{\DimSOSBasis}} )
{\Ite{\LyapMat}{0}}
=
\blockLyapMat
(\SOSmatBasis(\sysx) \otimes   {\Identity{\DimSOSBasis}} )
.
\end{align}	
Substituting this relation into \eqref{eq:SOSdotVpreMat_pd1} yields
\begin{align}
&
{\SOSdotVpreMat{\IDAffine}}(\sysx,{\Ite{\LyapMat}{0}},{\Ite{\InVecParam}{0}})
+
{\SOSdotVpreMat{\IDAffine}}(\sysx,{\Ite{\LyapMat}{0}},{\Ite{\InVecParam}{0}})^{\Mytop}
\nonumber\\&
= 
(\SOSmatBasis(\sysx) \otimes   {\Identity{\DimSOSBasis}} )^{\Mytop}
\blockLyapMat
{\solSOSRtdLyapMat{\IDAffine}}
\blockLyapMat
(\SOSmatBasis(\sysx) \otimes   {\Identity{\DimSOSBasis}} )
.\label{eq:SOSdotVpreMat_pd2}
\end{align}		
Because of the definition of $\SOSdotVBasis(\sysx)$, there exists a matrix $\connectMatF \in \{0,1\}^{\DimSOSmatBasis\DimSOSBasis \times \DimSOSdotVBasis } $ satisfying $\connectMatF \SOSdotVBasis(\sysx)=  	 \SOSmatBasis(\sysx) \otimes   \SOSBasis(\sysx)  $.
Combining this relation with  \eqref{eq:SOSdotVpreMat_pd2} yields
\begin{align}
&2 
\SOSBasis(\sysx)^{\Mytop}
{\SOSdotVpreMat{\IDAffine}}(\sysx,{\Ite{\LyapMat}{0}},{\Ite{\InVecParam}{0}})
\SOSBasis(\sysx)
\nonumber\\&
=
\SOSBasis(\sysx)^{\Mytop}
(\SOSmatBasis(\sysx) \otimes {\Identity{\DimSOSBasis}} )^{\Mytop}
\blockLyapMat
{\solSOSRtdLyapMat{\IDAffine}}
\blockLyapMat
(\SOSmatBasis(\sysx) \otimes {\Identity{\DimSOSBasis}} )
\SOSBasis(\sysx)
\nonumber\\&
=
(\SOSmatBasis(\sysx) \otimes\SOSBasis(\sysx) )^{\Mytop}
\blockLyapMat
{\solSOSRtdLyapMat{\IDAffine}}
\blockLyapMat 
(\SOSmatBasis(\sysx) \otimes \SOSBasis(\sysx) )
\nonumber\\&
=
\SOSdotVBasis(\sysx)^{\Mytop} 
\connectMatF^{\Mytop}
\blockLyapMat
{\solSOSRtdLyapMat{\IDAffine}}
\blockLyapMat
\connectMatF 
\SOSdotVBasis(\sysx)
.\label{eq:SOSdotVpreMat_pf}
\end{align}

Here, because $\sum_{\IDbEl=1}^{\DimSOSdotVBasis} {\El{  \connectMatF  }{\IDEl,\IDbEl}}=1$ holds for every $\IDEl$,
for every $\NotationVec \in \mathbb{R}^{\DimSOSdotVBasis}$, $\NotationVec \neq 0$ implies  $\connectMatF\NotationVec \neq 0$. 
Because ${\solSOSRtdLyapMat{\IDAffine}}\succ 0$ and $\blockLyapMat\succ 0$ hold,
for every $\NotationVec \in \mathbb{R}^{\DimSOSdotVBasis}$, $\NotationVec \neq 0$ implies
$
\NotationVec^{\Mytop}
\connectMatF^{\Mytop}
\blockLyapMat
{\solSOSRtdLyapMat{\IDAffine}}
\blockLyapMat
\connectMatF
\NotationVec
> 0
$.
This indicates 
$\connectMatF^{\Mytop}
\blockLyapMat
{\solSOSRtdLyapMat{\IDAffine}}
\blockLyapMat
\connectMatF \succ 0$.
Thus, using the following setting: 
\begin{align}
{\SOSdotVMat{\IDAffine}}
=
\connectMatF^{\Mytop}
\blockLyapMat
{\solSOSRtdLyapMat{\IDAffine}}
\blockLyapMat
\connectMatF /2
\succ 0
,\label{eq:def_dotLyapF_pd_sol}	
\end{align}
satisfies \eqref{eq:def_dotLyapF_new} and   \eqref{eq:def_symmetricSOSdotVMat} for 
$(\InVecParam, \LyapMat)=({\Ite{\InVecParam}{0}}, {\Ite{\LyapMat}{0}})$.
Meanwhile, for ${\SOSdotVMat{\IDAffine}}$ in \eqref{eq:def_dotLyapF_pd_sol}, there exists ${\ASeach{\IDAffine}}$ given by \eqref{eq:def_ASconstraint}.
By considering this ${\ASeach{\IDAffine}}$, Lemma \ref{thm:design_penaF} implies that the unique solution ${\SOSdotVMat{\IDAffine}}({\Ite{\InVecParam}{0}}, {\Ite{\LyapMat}{0}},\ASparam)$ to \eqref{eq:def_dotLyapF_new}--\eqref{eq:def_symmetricSOSdotVMat} is equal to  ${\SOSdotVMat{\IDAffine}}$  in \eqref{eq:def_dotLyapF_pd_sol} because of the uniqueness.
Therefore, there exist a feasible solution  $(\solASparam,{\dotVepsilon})$ to the SDP \eqref{eq:SDP_to_ASparam}.

Finally, we show  the statement {\MyLabelInitialVal}.	
The condition \eqref{eq:cond_robust_stab_P} is clearly satisfied by \eqref{eq:def_ini_LyapMat} because of $\solLyapInv\succ 0$.
The constraint in \eqref{eq:SDP_to_ASparam} indicates that a feasible solution  ${\Ite{\ASparam}{0}}=\solASparam$ satisfies ${\SOSdotVMat{\IDAffine}}({\Ite{\InVecParam}{0}}, {\Ite{\LyapMat}{0}},{\Ite{\ASparam}{0}})
\succ 0$ that is equivalent to \eqref{eq:cond_robust_stab_T}.
This completes the proof.

\input{TIsto_NLopt_Y_Ito_main_20250121d_submit_arxiv_bbl.bbl}
\end{document}

%% file: TIsto_NLopt_Y_Ito_symbols_20241223.tex

\newcommand{\SymColor}[1]{\textcolor{red}{#1}}  
\renewcommand{\SymColor}[1]{#1}  



\newcommand{\appcostLB}{\widehat{J}_{\mathrm{lb}}}


\newcommand{\MyLabelStoParamUnknown}{{(i)}}
\newcommand{\MyLabelPDFknown}{{(ii)}}


\newcommand{\MyLabelDifExplicit}{{(i)}}
\newcommand{\MyLabelDifMarkov}{{(ii)}}
\newcommand{\MyLabelDifStability}{{(iii)}}


\newcommand{\MyLabelFpoly}{{(i)}}
\newcommand{\MyLabelUpoly}{{(ii)}}
\newcommand{\MyLabeldotVpoly}{{(iii)}}
\newcommand{\MyLabelSOS}{{(iv)}}


\newcommand{\MyLabelFeasibleWP}{{(i)}}
\newcommand{\MyLabelFeasibleAddP}{{(ii)}}
\newcommand{\MyLabelInitialVal}{{(iii)}}



\newcommand{\PDF}[1]{\SymColor{p(}#1\SymColor{)}}
\newcommand{\MC}[2]{#1_{#2}}
\newcommand{\Ite}[2]{#1^{\SymColor{\{}#2\SymColor{\}}}}
\newcommand{\El}[2]{\SymColor{[}{#1}\SymColor{]}_{#2}}
\newcommand{\VEC}[2][]{\SymColor{\mathrm{vec}}#1({#2}#1)}

\newcommand{\VECH}[2][]{\SymColor{\mathrm{vech}}#1({#2}#1)}

\newcommand{\InvVEC}[4][]{\SymColor{\mathrm{vec}^{-1}_{#3 , #4}}#1({#2}#1)}
\newcommand{\TRACE}[1]{\SymColor{\mathrm{tr}}(#1)}
\newcommand{\Expect}[1]{\SymColor{\mathrm{E}}_{#1}}
\newcommand{\EigMin}[1]{\SymColor{\lambda_{\mathrm{min}}}(#1)}


\newcommand{\Mytop}{\SymColor{\top}}

\newcommand{\NotationVec}{\SymColor{\boldsymbol{v}}}
\newcommand{\NotationMat}{\SymColor{\boldsymbol{M}}}
\newcommand{\NotationSto}{\SymColor{\boldsymbol{\theta}}}
\newcommand{\NotationFunc}{\SymColor{\boldsymbol{\pi}}}
\newcommand{\IDNotation}{\SymColor{i}}
\newcommand{\IDbNotation}{\SymColor{j}}
\newcommand{\DimANotation}{\SymColor{n}}
\newcommand{\DimBNotation}{\SymColor{m}}

\newcommand{\Identity}[1]{\SymColor{\boldsymbol{I}}_{#1}}



\newcommand{\MyT}{\SymColor{t}}
\newcommand{\sysx}{\SymColor{\boldsymbol{x}}}
\newcommand{\IniState}{\SymColor{\boldsymbol{x}_{0}}}
\newcommand{\dotState}{\SymColor{\dot{\boldsymbol{x}}}}
\newcommand{\Input}{\SymColor{\boldsymbol{u}}}
\newcommand{\InpOpt}{\SymColor{\widehat{\Input}}}

\newcommand{\StoParam}{\SymColor{\boldsymbol{\theta}}}

\newcommand{\vertexStoParam}[1]{\SymColor{\boldsymbol{\theta}_{#1}}}

\newcommand{\Drift}{\SymColor{\boldsymbol{f}}}
\newcommand{\DriftMat}{\SymColor{\boldsymbol{F}}}
\newcommand{\InMat}{\SymColor{\boldsymbol{G}}}
\newcommand{\xDrift}[1]{\SymColor{\boldsymbol{f}_{#1}}}
\newcommand{\xInMat}[1]{\SymColor{\boldsymbol{G}_{#1}}}
\newcommand{\sDrift}[1]{\SymColor{h_{#1}}}
\newcommand{\PolyWeight}[1]{ \sDrift{#1}(\StoParam) }

\newcommand{\SetSto}{\SymColor{\mathbb{S}_{\StoParam}}}

\newcommand{\CoefReguVParam}{\SymColor{\eta}}
\newcommand{\ValReguVParam}{\SymColor{\widetilde{\VParam}}}

\newcommand{\PenaltyFuncRegu}{\SymColor{\rho}}

\newcommand{\costJ}{\SymColor{J}}

\newcommand{\MyDeltaFunc}{\SymColor{\delta_{\mathrm{D}}}}
\newcommand{\DeltaSampleState}[1]{\SymColor{\sysx^{(#1)}}}

\newcommand{\NumAffine}{\SymColor{K}}
\newcommand{\NumProj}{\SymColor{M}}
\newcommand{\NumMC}{\SymColor{S}}
\newcommand{\NumIteGrad}{\SymColor{N}_{\IDIteGrad}}

\newcommand{\IDProj}{\SymColor{m}}
\newcommand{\IDIteGrad}{\SymColor{\ell}}
\newcommand{\IDMC}{\SymColor{s}}

\newcommand{\IDEl}{\SymColor{i}}
\newcommand{\IDbEl}{\SymColor{j}}

\newcommand{\IDAffine}{\SymColor{k}}
\newcommand{\IDbAffine}{\SymColor{k^{\prime}}}

\newcommand{\CoefGrad}[1]{\SymColor{\alpha}_{#1}}

\newcommand{\InpParam}{\SymColor{\boldsymbol{W}}}
\newcommand{\OptInpParam}{\SymColor{\boldsymbol{W}_{\ast}}}
\newcommand{\VParam}{\SymColor{\boldsymbol{v}}}
\newcommand{\preVParam}{\SymColor{\widetilde{\boldsymbol{v}}}}

\newcommand{\EachValue}{\SymColor{\widehat{J}}}
\newcommand{\ErrV}{\SymColor{\epsilon_{J}}}

\newcommand{\uCost}{\SymColor{\boldsymbol{R}}}
\newcommand{\xCost}{\SymColor{q}}
\newcommand{\Bellman}{\SymColor{B}}

\newcommand{\xAterm}[1]{\SymColor{\boldsymbol{\psi}_{#1}(\sysx)}}
\newcommand{\xBterm}{\SymColor{\boldsymbol{\psi}_{0}(\sysx)}}

\newcommand{\MarginalMap}[2][]{\SymColor{\mathcal{M}#1(}#2\SymColor{#1)}}

\newcommand{\MarginalWeight}{\SymColor{\mathcal{P}}}

\newcommand{\NonArg}{\SymColor{\bullet}}

\newcommand{\InVecParam}{\SymColor{\boldsymbol{w}}}
\newcommand{\optInVecParam}{\SymColor{\boldsymbol{w}_{\ast}}}

\newcommand{\terminalT}{\SymColor{T}}


\newcommand{\LyapF}{\SymColor{V}}

\newcommand{\PolyDriftMat}[1]{\SymColor{\boldsymbol{F}}_{#1}}
\newcommand{\PolyInMat}[1]{\SymColor{\boldsymbol{G}}_{#1}}

\newcommand{\LyapMat}{\SymColor{\boldsymbol{P}}}
\newcommand{\LyapInv}{\SymColor{\boldsymbol{Q}}}
\newcommand{\solLyapInv}{\SymColor{\boldsymbol{Q}_{\ast}}}
\newcommand{\tdLyapMat}[1]{\SymColor{\boldsymbol{Y}_{#1}}}
\newcommand{\SOSRtdLyapMat}[1]{\SymColor{\boldsymbol{S}_{#1}}}
\newcommand{\solSOSRtdLyapMat}[1]{\SymColor{\boldsymbol{S}_{#1 \ast}}}
\newcommand{\SOSdotVMat}[1]{\SymColor{\boldsymbol{T}_{#1}}}

\newcommand{\tmpInpParam}{\SymColor{\boldsymbol{H}}}
\newcommand{\soltmpInpParam}{\SymColor{\boldsymbol{H}_{\ast}}}

\newcommand{\SOSAepsilon}{\SymColor{\epsilon}_{1}}
\newcommand{\solSOSAepsilon}{\SymColor{\epsilon_{1\ast}}}
\newcommand{\dotVepsilon}{\SymColor{{\epsilon}_{2}}}
%
%




\newcommand{\DEFSOSBasis}{\SymColor{\overline{\boldsymbol{z}}}}
\newcommand{\DEFSOSBBasis}{\SymColor{\overline{\boldsymbol{z}}^{\prime}}}
\newcommand{\DEFDimASOSBasis}{\SymColor{\overline{d_{1}}}}
\newcommand{\DEFDimBSOSBasis}{\SymColor{\overline{d_{2}}}}


\newcommand{\DimX}{\SymColor{d}_{\sysx}}
\newcommand{\DimU}{\SymColor{d}_{\Input}}
\newcommand{\DimSP}{\SymColor{d}_{\StoParam}}

\newcommand{\DimVParam}{\SymColor{d}_{\VParam}}
\newcommand{\DimInVecParam}{\SymColor{d}_{\InVecParam}}
\newcommand{\DimASparam}{\SymColor{d}_{\ASparam}}

\newcommand{\DimSOSBasis}{\SymColor{d}_{\SOSBasis}}
\newcommand{\DimRowInBasis}{\SymColor{d}_{\InBasis}}
\newcommand{\DimSOSdotVBasis}{\SymColor{d}_{\SOSdotVBasis}}
\newcommand{\DimSOSmatBasis}{\SymColor{d}_{\SOSmatBasis}}
\newcommand{\DimSOSdotVBBasis}{\SymColor{d}_{\SOSdotVBBasis}}

\newcommand{\SOSBasis}{\SymColor{\boldsymbol{z}}}
\newcommand{\InBasis}{\SymColor{\boldsymbol{Z}}}
\newcommand{\SOSmatBasis}{\SymColor{\boldsymbol{\zeta}}}
\newcommand{\SOSdotVBasis}{\SymColor{\boldsymbol{\xi}}}
\newcommand{\SOSdotVBBasis}{\SymColor{\widetilde{\boldsymbol{z}}}}

\newcommand{\uBasis}{\SymColor{\boldsymbol{\Phi}}}
\newcommand{\JBasis}{\SymColor{\boldsymbol{\phi}}}

\newcommand{\ASparam}{\SymColor{\boldsymbol{r}}}
\newcommand{\solASparam}{\SymColor{\boldsymbol{r}_{\ast}}}
\newcommand{\ASeach}[1]{\ASparam_{#1}}
\newcommand{\AScoefmat}{\SymColor{\boldsymbol{C}_{r}}}
\newcommand{\AScoefscalar}{\SymColor{{c}_{r}}}

\newcommand{\WolfeAcoef}{\SymColor{\chi}}
\newcommand{\ObjF}{\SymColor{g}}
\newcommand{\WolfeObjF}{\SymColor{\tilde{g}}}
\newcommand{\WolfeVar}{\SymColor{\boldsymbol{y}}}

\newcommand{\WolfeDiscount}{\SymColor{\gamma_{\alpha}}}
\newcommand{\WolfeIniStepSize}{\SymColor{\CoefGrad{\mathrm{ini}}}}

\newcommand{\simBasisDegA}{\SymColor{d_{1}}}
\newcommand{\simBasisDegB}{\SymColor{d_{2}}}

\newcommand{\degPow}{\SymColor{\boldsymbol{a}}}
\newcommand{\SOSdotVpreMat}[1]{\SymColor{\boldsymbol{U}}_{#1}}

\newcommand{\SOSpdMat}[1]{\SymColor{\overline{\boldsymbol{S}}}_{#1}}

\newcommand{\coeffdotVpreSqrtMat}[2]{\SymColor{\boldsymbol{c}}_{\IDAffine,#1,#2}}
\newcommand{\coeffVECdotVpreSqrtMat}{\SymColor{\boldsymbol{c}_{\IDAffine}(\InVecParam,\LyapMat)}}
\newcommand{\SOSalgebraicEq}[1]{\SymColor{\boldsymbol{b}}_{#1}}

\newcommand{\iniPenaltyFuncRegu}{\PenaltyFuncRegu_{\mathrm{ub}}}
\newcommand{\penaW}{\SymColor{\kappa}}


\newcommand{\BellResiCoef}{\SymColor{\beta}}
\newcommand{\BellA}{\SymColor{B_{\mathrm{qR}}}}
\newcommand{\BellB}{\SymColor{B_{\mathrm{J}}}}
\newcommand{\BellExact}{\SymColor{C_{\mathrm{\infty}}}}
\newcommand{\BellLim}{\SymColor{C_{\mathrm{T}}}}

\newcommand{\PROOFepsilon}{\SymColor{\varepsilon}}

\newcommand{\newDDmat}[3]{\SymColor{\boldsymbol{L}_{#2,#3}(#1)}}
\newcommand{\newDNmat}[2]{\SymColor{\boldsymbol{l}_{#2}(#1)}}

\newcommand{\StoDDmat}{\SymColor{\widetilde{\boldsymbol{L}}}}
\newcommand{\StoDNmat}{\SymColor{\widetilde{\boldsymbol{l}}}}


\newcommand*{\blockLyapInv}{\SymColor{\mathcal{Q}}}
\newcommand*{\blockLyapMat}{\SymColor{\mathcal{P}}}

\newcommand{\connectMatA}{\SymColor{\boldsymbol{M}_{\mathrm{a}}}}
\newcommand{\connectMatB}{\SymColor{\boldsymbol{M}_{\mathrm{b}}}}
\newcommand{\connectMatC}{\SymColor{\boldsymbol{M}_{\mathrm{c}}}}
\newcommand{\connectMatD}{\SymColor{\boldsymbol{M}_{\mathrm{d}}}}
\newcommand{\connectMatE}{\SymColor{\boldsymbol{M}_{\mathrm{e}}}}
\newcommand{\connectMatF}{\SymColor{\boldsymbol{M}_{\mathrm{f}}}}
\newcommand{\connectMatG}{\SymColor{\boldsymbol{M}_{\mathrm{g}}}}
